\documentclass[9pt]{article}
\usepackage{amsfonts,amssymb,amsmath}
\usepackage[T2A]{fontenc}  % in math
\usepackage[cp1251]{inputenc}
\usepackage[english]{babel}
\usepackage{geometry}
\usepackage{color}     %%%%%% ESEMPIO DI COME SCRIVERE A COLORI \textcolor{red}{paper}{paper}
\newtheorem{teor}{Theorem}[section]
\newtheorem{prop}[teor]{Proposition}
\newtheorem{lemma}[teor]{Lemma}

\numberwithin{equation}{section}

\usepackage{tocloft}

\def\ds{\displaystyle}
\def\f{\frac}
\def\div{{\rm div}}
\def\C{\mathcal{C}}
\def\R{\mathbb{R}}
\def\K{\mathcal{K}}
\def\L{\mathcal{L}}

\def\per{{\rm per}}

\def\n{\nabla}
\def\nab{\nabla}

\def\ve{\varepsilon}
\def\o{\omega}
\def\e{\varepsilon}
\def\C0{C^\infty_0(\R^d)}
\def\Cp{C_\per^\infty(\Box)}
\def\ld{L^2(\R^d)}
\def\rd{\R^d}
\def\RR^d{\R^d}
\def\ild{\int\limits_{\R^d}}
\def\ilb{\int\limits_{\Box}}
\def\sr{\lefteqn{~-}\int_}
\def\beq{\begin{equation}}
\def\eeq{\end{equation}}

\begin{document}
\title{On resolvent approximations   of elliptic
differential  operators with periodic coefficients  
}

\author{ S.\,E. Pastukhova}
%\begin{center}
%\begin{footnotesize} 
%Московский технологический университет
%
%Москва, Россия
%
%% email: pas-se@yandex.ru
%\end{footnotesize} 
%\end{center}
\date{}%не показывает
\maketitle

%\begin{abstract}

\begin{footnotesize} 
We consider  resolvents $(A_\e+1)^{-1}$
of %uniformly 
elliptic
second-order differential 
operators $A_\e=-\div\,a(x/\e)\nab$ in  $\R^d$ with $\e$-periodic %rapidly oscillating 
measurable matrix $a(x/\e)$
and study the asymptotic behaviour of $(A_\e+1)^{-1}$,
as the period $\e$ goes to zero. The matrix $a$ is % being 
not necessarily symmetric.
We provide a construction for the leading  terms of the “operator
asymptotics” of $(A_\e+1)^{-1}$
 in the sense of $L^2$-operator-norm convergence and prove
order $\e^2$ remainder estimates. We apply the %so-called
modified  method of the
first approximation with the usage of Steklov's smoothing.
The class of operators covered by our analysis
includes 
 uniformly elliptic families 
 with bounded coefficients
and also with  unbounded coefficients  from the John--Nirenberg space $BMO$ (bounded mean oscillation).

\end{footnotesize}
\bigskip

\section{Introduction}

\textbf{1.1. About the topic.}
This paper relates to  homogenization theory which studies heterogeneous media such as small-period composites or porous media in the limit of small period (for introduction to this theory see, for example, books \cite{BLP}--\cite{ZKO}). More precisely, the  paper relates to the rather new  branch  of homogenization theory connected with operator-type estimates for the error of homogenization. This topic 
attracts attention of many specialists last  decades; a lot of interesting results have been obtained through joint efforts of numerous mathematicians.

Among the pioneer publications devoted to operator-type estimates in homogenization of elliptic equations, we mention, first of all,  the papers
\cite{BS}--%, \cite{Zh1}, \cite{Zh05} and 
\cite{BS05}, where
a number of results have been established  concerning the difference, in
the operator $L^2$-norm, between the resolvent of the elliptic differential operator representing
the original heterogeneous medium depending on the small parameter $\e$, that is 
$$A_\ve =-\div\, a(x/\e)\nab,$$
and the resolvent of the operator
$$A_0 =-\div\, a^0\nab$$
 representing the %“homogenisation limit”
limiting  (or "effective"\,) medium, as $\e\to 0$. Here 
the matrix function $a$ is $[-1/2,1/2)^d$-periodic, symmetric, measurable, bounded and uniformly positive definite; the constant matrix
$a^0$ is of the same class, and it is found according a well known procedure. % as $a$.
 To study the difference between the resolvents
$(A_\e{+}1)^{-1}$ and $(A_0{+}1)^{-1}$ for the operators $A_\e$ and $A_0$ acting in the  space $L^2(\rd)$ means,
in other words, to study
the difference between the solutions to the elliptic problems
\beq\label{0.1}
u^\ve\in H^1(\RR^d),\quad A_\ve u^\ve+ u^\ve=f,\quad
f\in L^2(\RR^d),
\eeq
\beq\label{0.2}
u\in H^1(\RR^d),\quad A_0 u+ u=f,\quad
f\in L^2(\RR^d).
\eeq
%%%%%%%%%%%%%%%%%%%%%%
%%%%%%%%%%

 The uniform resolvent convergence of  $A_\ve$ to  $A_0$ in 
 $\ld$ %have been 
 was maintained, together with the rate of this convergence of order $\e$, in \cite{BS}, \cite{Zh1}. Thus, the resolvent 
$(A_0{+}1)^{-1}$ of the homogenized operator turns to be a good approximation for the resolvent
$(A_\e{+}1)^{-1}$ of the original operator in  $L^2$-operator norm with remainder term of order $\e$.
Naturally, the question arises about similar approximations of $(A_\e{+}1)^{-1}$
with remainder term of the next order, i.e., $\e^2$. More exactly, the question is what a correcting term of the form $\e\mathcal{C}_\e$ 
should be added to %the zero approximation
 $(A_0{+}1)^{-1}$ in order to attain the sharpness of order $\e^2$ for the approximation $(A_0{+}1)^{-1}+\e\mathcal{C}_\e$ 
 of $(A_\e{+}1)^{-1}$.
The answer on this question is also known, thanks to 
\cite{Zh05} and \cite{BS05}. The authors of both papers have
found such type %obtained  these 
approximations (in the framework of more general setups: including the case of  systems of elliptic equations in \cite{BS05} or the case of  equations in  $L^2$-spaces with  general Borel measures in \cite{Zh05})),
%with remainder term of order $\e^2$ 
acting by spectral method based on the Floquet--Bloch decomposition of the selfadjoint operator $A_\e$. Note that this approach is rather restrictive, for it  is closely linked  with  periodic problems 
since the Floquet--Bloch transformation works well exclusively in the case of operators with periodic coefficients. But 
homogenization theory is not limited only to periodic setup.

As in \cite{Zh05} and \cite{BS05}, \textit{
we analyse here the asymptotic behaviour of the resolvent
$(A_\e{+}1)^{-1}$ with the sharpness of order $\e^2$
in  $L^2$-operator norm, but 
under more general
conditions 
and by another method.} First, we allow the operator $A_\e$ to be nonselfadjoint with the matrix $a$ not necessarily symmetric which  entails more complicated structure of the correcting term  $\e\mathcal{C}_\e$ as compared with \cite{Zh05} and \cite{BS05}. Second, we relax the boundedness requirement in ellipticity condition on the matrix $a$ so that the approximation result remains the same  though additional arguments  are needed in justification of it. More precisely, the skew-symmetric part 
 of the diffusion matrix $a$ is allowed to be unbounded from the John--Nirenberg space $BMO$ (bounded mean oscillation).

Shortly, about the structure of the paper. The main results are formulated in 
%the end of \S3
theorems \ref{Th2.1}, \ref{Th2.2}  and \ref{Th6.1}. Their proof is given in \S5 and \S6. 
Sections \S\S1-3 are introductory, and \S4 and \S7 are devoted to the Steklov smoothing operator which plays the key role in our method.

\textbf{1.2. About the method.}
The present paper can be viewed as following in the footsteps of  \cite{Zh1} 
 in that it relies upon the so-called  "modified method of the first approximation" with the usage of the shift parameter (that is why it is called often shortly  as the shift method).
 %%%%
This method was proposed by V.V.Zhikov \cite{Zh1}
 as an alternative, along with the spectral approach used in in \cite{BS}, \cite{Zh05} and \cite{BS05}, to prove operator-type homogenization estimates; 
   it 
 turned to be universal in different setups:  periodic, locally periodic, quasiperiodic or multiscale.
The method has developed since 2005  in applications to various problems (we refer, e.g.,  to \cite{ZhP05}--\cite{PT17} and, in particular, to the overview \cite{UMN} where other  references are given). There have appeared two versions of the method: %this method: 
the original
 version with the usage of
the pure shift in the coefficients of the operator $A_\e$ (this creates a family of perturbated operators with a shift parameter $\o$, and %the subsequent 
averaging in $\o$ allows to overcome difficulties of estimating in the lack of the regularity for the data in the equation (\ref{0.1})),  and another  version with the usage of the Steklov smoothing operator (containing the shift implicitly as any other smoothing operator defined by means of convolution) 
%applied to the
embedded from the very beginning in the approximation sought.
 We use here the second version of the shift method. 
 
 Since 2005, when \cite{Zh05}, \cite{BS05}  and also \cite{Zh1} came up, it has been the
%mathematically 
challenge to obtain operator-norm resolvent-type
 homogenization estimates of order $\e^2$ 
from the
point of view close to the classical homogenization theory.
We recall that the error of homogenization for the equation (\ref{0.1})
is traditionally evaluated by means of direct  constructing approximations to the solution $u^\e$ via  two-scale expansions
\beq\label{0.3}
u^\e(x)\approx u^0(x, y)+\e u^1(x, y) + \e^2u^2(x, y) + \ldots , \quad y = x/\e,
\eeq
with functions $u^0(x, y)$, $u^1(x, y), \ldots$  periodic in $y$.
A regular way of finding %constructing
such functions is known. % was proposed. 
 For example, one should take the sum of three terms of the above two-scale expansion and try to enable
\beq\label{0.4}
(A_\e+1)(u^0+\e u^1+\e^2 u^2)-f=O(\e).
\eeq
%%%%
It is quickly seen that $u^0(x, y) = u(x)$ is
independent of $y$ and %this function %$u(x)$
 turns to be a solution to (\ref{0.2}). As for the next terms in the  two-scale expansion, we have
 $$u^1(x, y) = N^j(y)\f{\partial u(x) }{\partial x_j},\quad 
 u^2(x, y) = N^{ij}(y)\f{\partial^2 u(x) }{\partial x_i\partial x_j}
 $$
(summation over repeated indices is assumed from 1 to $d$). The
function $N^j$ here is the solution of the periodic problem on the  cell $%Y=
[-1/2,1/2)^d$  (see below (\ref{cp})). The function $N^{ij}$  is   the solution of another periodic problem on the  cell $%Y=
[-1/2,1/2)^d$ which we do not formulate in the present paper (see it, 
e.g., in Chapter IV of \cite{BP}).

The sum of the first two terms in the above expansion, namely,
\[
u^1_\e(x)= u(x)+\e N^j(x/\e)\f{\partial u(x) }{\partial x_j},
\]
is  usually  called the first approximation, $u(x)$ is the zero approximation,
and the term $\e N^j(x/\e)\f{\partial u(x) }{\partial x_j}$
is a corrector.
%%%%%%

All the  conclusions  derived here about the two-scale expansion (\ref{0.3})
are valid assuming that the matrix $a$ and  the right-hand side 
function $f$ are sufficiently regular. 
Under our minimal regularity conditions on the matrix $a$ and the function $f$, even the existence of $u^1_\e$ as an element of the space
$H^1(\RR^d)$ is under the question, and so inserting it into the original equation, as  in (\ref{0.4}), is impossible.

 Estimates of the form
\beq\label{0.5}
\ds{\|u^\e-u\|_{L^2}\le C\e, 
%\quad C=const(d,\lambda, \|f\|_{H^m}, \|a\|_{C^k}),
}\atop\ds{
\|u^\e-u^1_\e \|_{H^1}\le C\e}
\eeq
for the difference of  the solution $u^\e$ and its zero  and first approximations
were 
obtained long ago. 
However, the constant $C$ in such estimates appeared to depend on
the zero approximation $u$, which was %supposed 
to be sufficiently smooth. The latter is possible under relevant high regularity assumptions on the right-hand side function $f$.

Traditionally (see, e.g., Chapter  IV in \cite{BP}), the $H^1$-estimate (\ref{0.5})$_2$ was derived %proved 
at the first step from (\ref{0.4}) using %by means of 
the energy estimate
\[
\|v\|_{H^1}\le c\|(A_\e+1) v \|_{H^{-1}},\quad c=const(\lambda),
\]
and only then the $L^2$-estimate (\ref{0.5})$_1$ was deduced from (\ref{0.5})$_2$ as a simple corollary.
 Obviously, in this case the estimates (\ref{0.5}) do not admit operator formulation.

 Thus, to obtain the estimates (\ref{0.5}) 
 under our minimal regularity assumptions in  more or less
standard way, i.e., following in line with
 %a pattern of %classical
two-scale expansion method  described above, one should sufficiently modify %several items in 
the method. This was done in \cite{Zh1} and \cite{ZhP05} where two versions of the modified method of the first approximation were exposed for the first time.

\section{$L^2$- and $H^1$-estimates of order $\e$}
\textbf{2.1. $L^2$-esimates for the error of %convergence rate 
 homogenization.}
In the whole space
 $\R^d$, $d\ge 2$, consider a  
 divergent-type second order elliptic equation 
\begin{equation}\label{1}\ds
u_\ve \in H^1 (\R^d), \quad 
{A_{\ve}u^\ve+u^\ve=f,\quad f{\in} L^2(\R^d),}
\atop\ds{A_{\ve}=-\div\,a_\ve(x)\nab,\quad a_\ve(x)=a(\ve^{-1}x),}
\end{equation}
with a small parameter  $\ve\in(0,1)$. Coefficients of the equation are $\ve$-periodic
and, thus, are rapidly oscillating as $\ve\to 0.$
Here
$a(x){=}\{a_{jk}(x)\}_{j,k=1}^d$ 
  is a measurable 
 $1$-periodic matrix with real entries. The periodicity cell is the unit cube 
  $\Box=[-\f{1}{2},\f{1}{2})^d$.
We suppose that
\beq\label{2}
\lambda |\xi|^2\le {a}\xi\cdot\xi
,\quad {a}\xi\cdot\eta
\le\lambda^{-1}|\xi|\,|\eta|
\quad\forall \xi,\eta\in\R^d
\eeq
for some $\lambda>0$. The matrix $a$ is not necessarily symmetric.

We associate with 
(\ref{1}) 
the homogenized equation
\beq\label{3}\ds{
u \in H^1 (\R^d), \quad 
 A_0 u+u
 =f,
 }\atop\ds{
 A_0=-\div \,{a^0}\nab,
 }
\eeq
where ${a^0}$ is a constant matrix of the same class 
  (\ref{2}); %and it 
  ${a^0}$ is calculated %can be found  
  according to the well known procedure in terms of solutions to  auxiliary
periodic problems (see below (\ref{cp}), (\ref{13})).
 
 Solutions to (\ref{1}) and (\ref{3}) are understood in the sense of distributions in
  $\R^d$.
For example, as for (\ref{1}),
the following integral identity holds
\beq\label{4}%\[
(A_{\ve}u^\ve+u^\ve,\varphi)=\ild (a_\ve(x)
\nab u^\ve\cdot \nab\varphi+u^\ve\varphi)\,dx=\ild f\varphi\,dx\quad \forall\varphi\in \C0.
\eeq%]
By the closure, the test functions here can be taken from the space
  $H^1 (\R^d)$. In particular, inserting %the solution
   $\varphi=u^\ve$ in (\ref{4}) yields the energy inequality
\begin{equation}\label{5}\ds{
\lambda\|\n u^\ve\|^2+\| u^\ve\|^2\stackrel{(\ref{2})_1}\le
(a_\ve \n u^\ve ,\n u^\ve)+(u^\ve ,u^\ve)= (f,u^\ve)\le
\| f\| \| u^\ve\|,
} \atop\ds{
\|u^\ve\|\le \|f\|,\quad \lambda \|\nab u^\ve\|^2\le \|f\|^2.
}
\end{equation}
Here and in what follows,  we use the simplified 
notation for the inner product and the norm  in $\ld$
 %\[%
 \begin{equation}\label{6}
\|\cdot\|=\|\cdot \|_{\ld},\quad
(\cdot \,,\cdot\,)=(\cdot \,,\cdot\,)_{\ld}.
\end{equation}%\]

The unique solvability of the equation is established by the Lax–Milgram lemma. This fact is true   actually  %in fact
 for more general right-hand side functions, namely, for $f\in  H^{-1} (\R^d)$ (where $ H^{-1} (\R^d)$ is the dual of 
$ H^1 (\R^d)$), so that the resolvent $(A_\e{+}1)^{-1}: H^{-1} (\R^d)\to  H^1 (\R^d)$ %exists %and is
is a bounded operator. The same is valid for
the resolvent $(A_0{+}1)^{-1}$. But if its action is restricted on the space
$L^2(\R^d)$, the property $(A_0{+}1)^{-1}: L^2(\R^d)\to  H^2 (\R^d)$ is gained.
In other words,
the elliptic estimate holds
 for the solution to the homogenized equation:
\beq\label{18}
\|u \|_{ H^2(\R^d)}\le c \|f\|, \quad 
c=cost(\lambda), 
\eeq
which can be easily established by means of the Fourier transform
because the matrix $a^0$ is constant and positive definite.

The homogenization result for (\ref{1}) is known from long ago and can be formulated, for example, as
$G$-convergence of operators $A_\e$ to $A_0$ (see \cite{ZKOH} and references therein) which means that
\beq\label{18a}%\[
\lim_{\e\to 0}\langle h,(A_\e{+}1)^{-1}f\rangle=\langle h,(A_0{+}1)^{-1}f\rangle
\eeq%\]
for any $f,h\in H^{-1} (\R^d)$.   Here the value of a functional $h\in
 H^{-1} (\R^d)$ at a $v\in  H^1 (\R^d)$ is denoted by $\langle h,v\rangle$.  In other words, (\ref{18a}) means that,  for an arbitrary right-hand side function $f\in H^{-1} (\R^d)$, the solutions of equations (\ref{1}) and (\ref{3}) are connected with the the  weak  convergence in $ H^1 (\R^d)$ and, as a corollary, with the the  weak  convergence in $L^2 (\R^d)$. 
 From here by the energy method and lower semicontinuity arguments, one can  derive the strong convergence
$u_\e\to u$ in $ \ld$ which means in operator terms 
%implies % in other words, 
the strong resolvent convergence 
$$(A_\e{+}1)^{-1}\to (A_0{+}1)^{-1}\quad \mbox { in } 
\quad \ld .$$
This operator convergence  can be further strengthened up to the uniform resolvent convergence
with the following rate convergence estimate %
\beq\label{8}
\|(A_\e+1)^{-1}-(A_0+1)^{-1}\|_{\ld\to \ld}\le c\ve,\quad c=const(d,\lambda).\eeq
%which 
One can rewrite (\ref{8}) in terms of the solutions to (\ref{1}) and (\ref{3}) as follows
\beq\label{9}
\|u^\e-u\|\le c\ve \|f\| 
\eeq
with the same  right-hand side constant $c$ depending only on the dimension $d$ and the ellipticity constant $\lambda$ from (\ref{2}).
To prove the estimate (\ref{8}) in the self-adjoint case the authors of \cite{BS} 
used the spectral approach based on operator-theoretic arguments tightly bound to the self-adjoint situation.         %To this end
Quite different method to prove (\ref{8}) % for this goal
 was proposed, first, 
in \cite{Zh1} and then developed in \cite{ZhP05}. This is the modified method of the first approximation 
 with the usage of shift or smoothing
operators. From the very beginning of the appearance, this method turned out  to be  universal for studying 
various homogenization problems 
which admit 
nonselfadjointness, nonlinearity, divergence-form and nondivergence-form equations, different types of degeneracy, high order or vector equations, and others 
(see, e.g. \cite{Zh1}-\cite{PT17} and  also references in the overview \cite{UMN}). 

\textbf{2.2. Homogenization attributes.} Consider the following periodic problem on the unit cube $\Box=[-\f{1}{2},\f{1}{2})^d$
\beq\label{cp}\ds{
N^j\in H_{\rm per}^1(\Box),\quad\div_y a(y)(e^j+\nab_y N^j)=0, }\atop\ds
{\langle N^j\rangle=0,\quad 
j=1,  ..., d,}
\eeq
where  $e^1, \dots , e^d$ is a canonical basis in $\R^d$,
$H_{\rm per}^1(\Box)$ is the Sobolev space of 1-periodic functions,
% and
 $$\langle \cdot\rangle=\int\limits_\Box\,\cdot\, dy.$$ 
 Then the homogenized matrix $a^0$ is defined
 in terms of the solutions to the cell problem (\ref{cp}) by equalities
\beq\label{13}
{a^0e^j}=\langle {a}(e^j+\nab N^j)\rangle,\quad j=1, \dots ,d.
\eeq

A solution to the problem (\ref{cp}) is understood in the sense of the integral identity for smooth periodic functions
\beq\label{14}%$$
\langle a(e^j+\nab N^j)\cdot\nab \varphi\rangle=0,\quad \varphi\in \Cp,
\eeq%$$
which can be extended by closure to test functions in 
$H_{\rm per}^1(\Box)$. On the other hand, Equation (\ref{cp}) can be regarded in the sense of distributions in
  $\R^d$, that is a known fact in homogenization theory. Thus, together with
(\ref{14})  the integral identity is satisfied with test functions
in   $\C0$.

Introduce the  1-periodic vector 
 \begin{equation}\label{g}
g^j(y):=a(y)\left(\nab N^j(y)+e^j\right)-a^0 e^j,\quad j=1,\ldots,d.
\end{equation}
which is solenoidal and has zero mean value, i.e.,
 \begin{equation}\label{sol}
\div\, g^j(y)=0,\quad \langle g^j\rangle=0,
\end{equation} 
by  (\ref{14}) and (\ref{13}) respectively.
The property (\ref{sol})$_1$ may be understood in both ways: in the sense of the integral identity of the type (\ref{14}) or in the sense of distributions in
$\R^d$.
 
 Let $A_\e^*$ be the adjoint of
  $A_\e$ and consider the problem
 \begin{equation}\label{1s}\ds
v^\ve \in H^1 (\R^d), \quad 
{A_{\ve}^*v^\ve+v^\ve=h,\quad h{\in} L^2(\R^d),}
\atop\ds{A_{\ve}^*=-\div\,a^*_\ve(x)\nab,\quad a^*_\ve(x)=a^*(\ve^{-1}x),}
\end{equation}
where 
$a^*$ is the transpose of $a$.

It is known that the homogenized equation for
(\ref{1s}) 
 will be  
 \beq\label{3s}
v \in H^1 (\R^d), \quad 
 A^*_0 v+v=-\div \,{(a^0)^*}\nab v+v=h,
\eeq
where $A^*_0$ is the adjoint of $A_0$ and has the matrix   $(a^0)^*$ transposed to ${a^0}$.
Thus,
\beq\label{hom}
(a^*)^0 =(a^0)^*.
\eeq

 The counterpart of  the cell problem  (\ref{cp}) will be
\beq\label{cps}\ds{
\tilde{N} ^j\in H_{\rm per}^1(\Box),\quad\div_y a^*(y)(e^j+\nab_y \tilde{N}^j)=0, }\atop\ds
{\langle \tilde{N}^j\rangle=0,\quad 
j=1,  ..., d.}
\eeq
Its solutions  generate formally  the homogenized matrix for the equation (\ref{1s}) 
 through the formula similar to (\ref{13}), and so $\tilde{N} ^j$  are connected with the matrix ${a^0}$, namely, 
\beq\label{13s}
{(a^0)^*e^j}=\langle a^*(e^j+\nab \tilde{N}^j)\rangle,\quad j=1, \dots ,d,
\eeq
where (\ref{hom}) is taken into account.

We introduce also  the counterpart of 
(\ref{g}) 
 \begin{equation}\label{gs}
\tilde{g} ^j(y):=a^*(y)\left(\nab \tilde{N}^j(y)+e^j\right)-(a^0)^* e^j,\quad j=1,\ldots,d,
\end{equation}
which satisfies the relations  
  \begin{equation}\label{sols}
\div\, \tilde{g}^j(y)=0,\quad \langle \tilde{g}^j\rangle=0,
\end{equation}
by (\ref{cps}) and (\ref{13s}). 

In the sequel, we will refer to the energy and elliptic estimates relating to (\ref{1s}) and (\ref{3s}) respectively, those are
\beq\label{5s}
\|v^\ve \|_{ H^1(\R^d)}\le c \|f\|, \quad 
c=cost(\lambda), 
\eeq
\beq\label{18s}
\|v \|_{ H^2(\R^d)}\le c \|f\|, \quad 
c=cost(\lambda). 
\eeq
 
\textbf{2.3. $H^1$-approximations in homogenization.} According to (\ref{8}),
the resolvent $(A_0+1)^{-1}$ approximates $(A_\e+1)^{-1}$ in 
$L^2$-operator norm with the error of order $\e$.
If the resolvent
 $(A_\e+1)^{-1}$ is regarded as an operator from $\ld$
to $ H^{1} (\R^d)$, then  for  its approximation we need the sum
$(A_0+1)^{-1}+\e \K_\ve$, where $\K_\ve$ is a correcting operator, and so
\beq\label{10}
\|(A_\e+1)^{-1}-(A_0+1)^{-1}-\e \K_\ve\|_{\ld\to H^{1} (\R^d)}\le c\ve,\quad c=const(d,\lambda). 
\eeq
 The correcting operator $\K_\ve:\ld\to H^{1} (\R^d)$ is defined by
 \beq\label{2.20}
\K_\ve f = N_\e\cdot\nab S^\e (A_0+1)^{-1}f,
%\quad N_\e(x)=N(\e^{-1}x),
%(\frac{\cdot}{\e})
\eeq
where $ N_\e(x)=N(\e^{-1}x)$, $N(y)=\{N^j(y)\}_{j=1}^d$ %,\ldots,N^d(y))$ 
is the periodic vector composed of the solutions to (\ref{cp}) and 
$S^\e$ is the Steklov smoothing operator (see the definition of  $S^\e$ in \S 4 below).
Then 
\[
\|\e\K_\ve\|_{\ld\to H^{1} (\R^d)}\le c, \quad 
c=cost(d,\lambda),
\]
in view of properties  of the smoothing operator (see Lemma \ref{LemM1}) and the elliptic estimate (\ref{18}).

In the scalar case, the correcting operator can be constructed
 without  smoothing. 
Letting
\beq\label{111}
K_\ve f= N_\e\cdot\nab (A_0+1)^{-1}f, 
\eeq
we have, instead of (\ref{2.20}),  the  operator $ K_\e:\ld\to H^1(\R^d)$  such that
\beq\label{17}
\|\e K_\e f\|_{H^1(\R^d)}\le c\|f\|,\quad c=cost(d,\lambda).
\eeq
Estimate (\ref{17}) implies that the norms
$\|\e N_\e\cdot\nab u\|$ and
$\|\nab(\e N_\e\cdot\nab u)\|$, where $u$ is the solution to the homogenized equation,
 are finite  and $\varepsilon$-uniformly bounded by
$\|f\|$. This fact is not at all obvious, but it takes place because the solution to the cell problem belongs actually to the space 
$L^\infty(\Box)$ 
in view of the generalized maximum principle which is valid  for scalar equations, but not for vector ones. What is more,  the boundedness of the solution
$N^j$  entails the %following
 multiplier property of its gradient
$$\nab N^j:  H^1(\R^d)\to \ld$$ with the estimate
\beq\label{11m}%\[
\|(\nab N^j)(x/\e)z\|^2_{\ld}\le C\left(\|z\|^2_{\ld}+\e^2\|\nab z\|^2_{\ld}
\right),\quad z\in H^1(\R^d),
\eeq%\]
where the constant $C$ depends only on 
the constant $\lambda$ in (\ref{2}).  As a result, we have the boundedness property
(\ref{17}) and the following estimate holds
\beq\label{11}
\|(A_\e+1)^{-1}-(A_0+1)^{-1}-\e K_\ve\|_{\ld\to H^{1} (\R^d)}\le c\ve,\quad c=const(d,\lambda), 
\eeq
with the correcting operator $K_\ve$ defined in (\ref{111}).

The operator estimates (\ref{10}) and (\ref{11}) were first proved 
in \cite{Zh2} and \cite{ZhP05} by usage of  shifting or smoothing 
respectively.

Since the smoothing operator $S^\e$ is included in the corrector, it is possible not only well define 
$ H^1$-approximation, but also to overcome technical difficulties to estimate its residual in the equation. These difficulties arise under the minimal regularity  conditions on the data of the problem (\ref{1}). Here, we essentially use the properties of the smoothing operator  
   $S^\e$  relating to  its interaction
  with $\e$-periodic factors (see \S 4). These properties were first  
  noticed in \cite{ZhP05},  \cite{ZhP06}.

\section{$L^2$-estimate of order $\e^2$}

The operator $\K_\e$ defined in (\ref{2.20}) is a bounded operator in
$\ld$ with the estimate for the norm %, moreover,
 $\|\K_\ve\|_{\ld\to \ld}\le c$, and  the adjoint operator
$(\K_\ve)^*:  \ld \to \ld$ is such that
 \begin{equation}\label{2.0021}%$$
%\e N_\e\cdot S^\e  \nab (A_0+1)^{-1}f=:
(\K_\ve)^*f:= (A_0+1)^{-1} S^\e \div( N_\e\, f).
\end{equation}%$$

%In the self-adjoint case, where 
Suppose that the matrix $a$ is symmetric. Then the sum
$\e\K_\ve+\e(\K_\ve)^*$ turns to be the true
%is a  в $ \ld$, имеет норму порядка $\e$ 
 correcting operator of $(A_0+1)^{-1}$ in approximations
with remainder of order $\e^2$ for the resolvent $(A_\e+1)^{-1}$ in   
$L^2$-operator norm. 
The following   estimate holds:
 \begin{equation}\label{2.21}\ds{
\|(A_\e+1)^{-1}-(A_0+1)^{-1}-\e\K_\e -\e(\K_\e)^* \|_{\ld\to \ld}\le C\ve^2, }
\atop\ds{
\K_\e= N_\e\cdot S^\e  \nab (A_0+1)^{-1},}
\end{equation}
where the constant $C$ depends only on the dimension $d$ and the ellipticity constant $ \lambda$.

Since, under the assumption (\ref{2}) in the scalar case, the solution $N^j$ to the problem (\ref{cp}) belongs to $L^\infty(\Box)$ in view of the generalized maximum principle,  
%заданный в (\ref{2.20}) 
the operator $\K_\e$ in the estimate (\ref{2.21})  can be replaced with the simpler operator $K_\e$ defined in
(\ref{111}). Thus, the following   estimate holds:
 \begin{equation}\label{2.22}\ds{
\|(A_\e+1)^{-1}-(A_0+1)^{-1}-\e K_\e -\e(K_\e)^* \|_{\ld\to \ld}\le C\ve^2,}\atop\ds{ 
K_\e= N_\e\cdot   \nab (A_0+1)^{-1},}
\end{equation}
with the constant $C$ of the same type as in (\ref{2.21}).

 The estimate (\ref{2.21}) was proved in \cite{P20}
 by using the modified
  method of the first approximation,  and the estimate (\ref{2.22}) was derived from (\ref{2.21}) as a simple corollary by properties of smoothing.
We make some remarks on these estimates.

%\noindent
%\textbf{Замечания.}

1)
Pay attention  on 
the selfadjointness of the  both 
approximations  (\ref{2.21}) and (\ref{2.22}) which is 
% turn to be selfadjoit operators on 
contrary to the situation in  (\ref{10}) and (\ref{11}), where the correcting operators are not selfadjoint.

2) The estimate quite similar to (\ref{2.22}) was proved within the framework of more general results in  \cite{Zh05} and \cite{BS05} by using the spectral approach based on the Bloch--Floquet decomposition of selfadjoint differential operators with periodic coefficients.

 3)  The estimate resembling (\ref{2.21}), but with the smoothing operator $\Pi^\e$ of  another type, was obtained both in \cite{Zh05} and \cite{BS05}. 
 The pseudodifferental operator $\Pi^\e$ acting as 
  \begin{equation}\label{Pro}
  \Pi^\e\varphi(x)=F^{-1}\left(
 1_{\{|\xi|\le 1/\e\}}(F\,\varphi)(\xi)\right)
 \end{equation}
  naturally arises
 within the scope of the spectral method.
 Here $F$ denotes the Fourier transform and  $F^{-1}$ is its inverse,
$1_{\{|\xi|\le 1/\e\}}$ is a characteristic function  of the cube
$\{\xi:\,|\xi|\le 1/\e\}$.  Evidently, $\Pi^\e$ has smoothing properties, though it emerges as a result of some projection.

\bigskip
 Suppose now that
   the matrix $a(y)$ in (\ref{1}) is not symmetric. %
Then the correcting operator in  approximations of the resolvent  $(A_\e+1)^{-1}$ with  remainder  of order $\e^2$ will be more complicated than in (\ref{2.21}) and it is constructed of three terms: one of them does not contain oscillating factors, and the remaining two terms are similar to those in (\ref{2.21}).

 \begin{teor}\label{Th2.1} Let 
 $N(y)=\{N^j(y)\}_{j=1}^d$, %,\ldots,N^d(y))$ 
$\tilde{N} (y)=\{\tilde{N}^j(y)\}_{j=1}^d$ be the vectors composed of solutions to (\ref{cp}) and (\ref{cps}), and 
$S^\e$ be the Steklov smoothing operator (see (\ref{m.0})).
 Then the following   estimate holds for the resolvents 
 %of the operators 
 $(A_\e+1)^{-1}$ and $(A_0+1)^{-1}$ of the problems (\ref{1}) and (\ref{3}):
 \begin{equation}\label{2.021}%\ds{
\|(A_\e+1)^{-1}-(A_0+1)^{-1}-\e\K_\e -\e(\tilde{\K}_\e)^* -\e\L\|_{\ld\to \ld}\le C\ve^2, 
%} \atop\ds{ }
\end{equation}
where
 \begin{equation}\label{2.121}\ds{
\K_\e= N(\f{.}{\e})\cdot S^\e  \nab (A_0+1)^{-1},\quad 
\tilde{\K}_\e=\tilde{N}(\f{.}{\e})\cdot S^\e  \nab (A^*_0+1)^{-1}
, }
\atop\ds{\L= (A_0+1)^{-1}\left( \tilde{c}_i^{jk}-c_i^{jk}\right)\f{\partial^3 }{\partial x_j\partial x_i\partial x_k}
%(\mathcal{C}
%+\tilde{\mathcal{C}}^*) 
(A_0+1)^{-1} %\quad 
%\mathcal{C}
%+\tilde{\mathcal{C}}^*=
%\left( \tilde{c}_i^{jk}-c_i^{jk}\right)\f{\partial^3 }{\partial x_j\partial x_i\partial x_k}
}
\end{equation}
and the constant coefficients $\tilde{c}_i^{jk}$, $c_i^{jk}$ are defined in (\ref{3.24})
in terms of the functions
$N^j$, $\tilde{N}^j$ and its gradients.\\
${}$ ${}{}$ The constant $C$ in (\ref{2.021}) depends only on the dimension $d$ and the ellipticity constant $ \lambda$.\\
${}$ ${}{}$ %In the %selfadjoint case
If the matrix $a(y)$ is symmetric, the approximation for $(A_\e+1)^{-1}$ defined in (\ref{2.021}) and (\ref{2.121}) 
reduces into that of (\ref{2.21}).
\end{teor}

In the scalar case, the solutions $N^j$ and $\tilde{N}^j$ to the cell problems  belong to $L^\infty(\Box)$ in view of the generalized maximum principle, and so the smoothing operator in
the approximation from (\ref{2.021}) can be dropped.

\begin{teor}\label{Th2.2} 
 The estimate (\ref{2.021}) remains valid if the smoothing operator
 $S^\e$ is omitted in correcting operators (\ref{2.121}). 
 \end{teor}

 Theorems  \ref{Th2.1} and \ref{Th2.2} are proved in \S 5. In \S4 we introduce the Steklov smoothing operator and list its properties  that are applied in our considerations. Some of these  properties
 have not been noticed before, and so they are proved in  \S7.
 
 \textbf{Remark 3.3.} The  results similar to theorems \ref{Th2.1} and \ref{Th2.2} are proved in \cite{Se1} with the difference that, instead of the Steklov  smoothing operator $S^\e$, the smoothing operator (\ref{Pro}) is embedded in the correcting terms (\ref{2.121}). 
 The  operator (\ref{Pro}) appears there in the corrector 
 just like in \cite{Zh05} and \cite{BS05} as a by-product of applying the  Floquet--Bloch transformation 
 with the purpose to reduce the problem in the whole space  $\rd$ to the problem on the cell of periodicity $\Box=[-1/2,1/2)^d$.

\textbf{Remark 3.4.}  It is  worth noting that, once the  estimate (\ref{10}) 
   in the operator $(L^2\to H^1)$-norm with order $\e$ remainder is verified, 
 the estimate of the type (\ref{2.021}) 
 (or,  in selfadjoint setup,    its variant (\ref{2.21}) with the simpler corrector)
 in the operator $(L^2\to L^2)$-norm with order $\e^2$ remainder is surely guaranteed by the method we demonstrate here.

\textbf{Remark 3.5.}
 In the present paper, we restrict ourselves to the scalar case
 only for the sake of simplicity. 
 We deal with %consider
  the classical
diffusion equation of the type (\ref{0.1}) or its appropriate perturbations.
 Although  the maximum principle
is  valid in the scalar case, it is not used in our constructions and in the main proof, and so
 the result also carries over to vector models, including, e.g., the elasticity theory system or other systems considered in \cite{P17}.
 
\section{Properties of the smoothing operator}

In our method, the Steklov smoothing operator
 \begin{equation}\label{m.0}% $$
S^\e\varphi(x)=\ilb\varphi(x-\e\omega)\,d\omega 
\end{equation}% $$
called also the Steklov average,
plays the key %central 
role, as it was already explained in \S 2 and \S3.
We begin with the simplest and most known properties of this operator:
  \begin{equation}\label{m.1}
\|S^\e\varphi\|_{L^2(\rd)}\le\|\varphi\|_{L^2(\rd)},
\end{equation}
 \begin{equation}\label{m.2}
\|S^\e\varphi-\varphi\|_{L^2(\rd)}\le (\sqrt{d}/2)\e\|\nab\varphi\|_{L^2(\rd)},
\end{equation}
 \begin{equation}\label{m.3}
\|S^\e\varphi-\varphi\|_{H^{-1}(\rd)}\le(\sqrt{d}/2)\e\|\varphi\|_{L^2(\rd)}.
\end{equation}
To supplement (\ref{m.1}) note that  $S^\e$
is a selfadjoint operator in $L^2(\rd)$.
We also mention the obvious property $S^\e(\nab \varphi)=\nab 
(S^\e\varphi)$, thereby, $S^\e$ and any differential operator with constant coefficients commute with one another. As a corollary,  $S^\e$ commutes with the resolvent $(A_0+1)^{-1}$ either.

%Besides, 
 The following properties of the Steklov smoothing are displayed
 in interaction with 
 $\e$-periodic factors.
 \begin{lemma}\label{LemM1} If $\varphi{\in}L^2(\R^d)$, $b{\in}L^2_\per(\Box)$,
  $b_\e(x){=}b(\e^{-1}x)$,
 then $b_\e S^\e\varphi\in L^2(\R^d)$ and
 \begin{equation}\label{m.4}
\|b_\e S^\e\varphi\|^2_{L^2(\rd)}\le\langle b^2\rangle\|\varphi\|^2_{L^2(\rd)}.
\end{equation}
\end{lemma}

 \begin{lemma}\label{LemM2} If $b{\in}L^2_\per(\Box)$, $\langle b\rangle{=}0$, $b_\e(x){=}b(\e^{-1}x)$,
  $\varphi{\in}L^2(\rd)$,
  $\Phi{\in}H^1(\rd)$,
 then
 \begin{equation}\label{m.5}
( b_\e S^\e\varphi,\Phi) \le C\e
\langle b^2\rangle^{1/2}\|\varphi\|_{L^2(\rd)}\|\nab\Phi\|_{L^2(\rd)},\quad C=const(d).
\end{equation}
\end{lemma}

 The properties (\ref{m.4}), (\ref{m.5}) were 
 %discovered  noticed first 
 highlighted and proved in \cite{ZhP05}, \cite{ZhP06} (see also \cite{UMN}).
 
We formulate the assertions of Lemmas \ref{LemM1} and \ref{LemM2} in the operator form.
 \begin{lemma}\label{LemM3}
Under the conditions of Lemma \ref{LemM1}, the norms of the operators $b_\e S^\e:L^2(\rd)\to L^2(\rd)$ are uniformly bounded:
 \begin{equation}\label{m.40}
\|b_\e S^\e\|_{L^2(\rd)\to L^2(\rd)}\le \langle b^2 \rangle^{1/2}.
\end{equation}
Furthermore, if $\langle b\rangle{=}0$, then $b_\e S^\e:L^2(\rd)\to H^{-1}(\rd)$ and
 \begin{equation}\label{m.50}
\|b_\e S^\e\|_{L^2(\rd)\to H^{-1}(\rd)}\le C\varepsilon\langle b^2 \rangle^{1/2},\quad C=const(d).
\end{equation}
\end{lemma}
 
 The estimates (\ref{m.2}) and (\ref{m.5}) can be specified under  assumptions of higher regularity.
 
For example, if $\varphi\in H^2(\R^d)$, then
 \begin{equation}\label{m.7}
\|S^\e\varphi-\varphi\|_{L^2(\rd)}\le C\e^2\|\nab^2\varphi\|_{L^2(\rd)},\quad C=const(d).
\end{equation}
Indeed, we write the equality
\[
\varphi(x+h)-\varphi(x)-\nab \varphi(x)\cdot h=
\int\limits_{0}^1(1-t)\nab(\nab\varphi(x+th)\cdot h)\cdot h\,dt
\]
and, setting  $h=-\e\omega$, integrate it over  
$\omega\in\Box=[-\f{1}{2},\f{1}{2})^d$. As a result, we arrive at the integral representation for the difference $S^\e\varphi-\varphi$ in terms of the second order gradient $\nab^2\varphi$.
Consequently,  
\[
|S^\e\varphi(x)-\varphi(x)|\le
\varepsilon^2\ilb
\int\limits_{0}^1|\nab(\nab\varphi(x-t\e\omega)\cdot\omega)\cdot\omega|\,dt\,d\omega,
\]
which imlies (\ref{m.7})  by the  H\"{o}lder inequality.

As for Lemma \ref{LemM2}, its extension will be
\begin{lemma}\label{LemM4} Assume that $b{\in}L^2_\per(\Box)$, $\langle b\rangle{=}0$, $b_\e(x){=}b(\e^{-1}x)$,
  $\varphi,\psi{\in}H^1(\rd)$. Then
  \begin{equation}\label{m.8}
( b_\e S^\e\varphi,S^\e\psi) \le C\e^2
\langle b^2\rangle^{1/2}\|\nab\varphi\|_{L^2(\rd)}\|\nab\psi\|_{L^2(\rd)},\quad C=const(d).
\end{equation}
\end{lemma}

The further extension is given  by
%Lemma \ref{LemM4} can be generalized as follows.
\begin{lemma}\label{LemM5} Assume that $\alpha,\beta\in L^2_\per(\Box)$, $\langle \alpha\beta\rangle=0$, $\alpha_\e(x)=\alpha(\e^{-1}x)$, $\beta_\e(x)=\beta(\e^{-1}x)$,
  $\varphi,\psi\in H^1(\rd)$. Then
  \begin{equation}\label{m.13}
( \alpha_\e S^\e\varphi,\beta_\e S^\e\psi) \le C\e^2
\langle \alpha^2\rangle^{1/2}\langle \beta^2\rangle^{1/2}
\|\nab\varphi\|_{L^2(\rd)}\|\nab\psi\|_{L^2(\rd)},\quad C=const(d).
\end{equation}
\end{lemma}

 Note that the form $( \alpha_\e S^\e\varphi,\beta_\e S^\e\psi)$ in (\ref{m.13}) is well defined since both functions $\alpha_\e S^\e\varphi$ and $\beta_\e S^\e\psi$ belong to $ L^2(\rd)$, by Lemma  \ref{LemM1}.

%We formulate also the following 
Another  extension of Lemma \ref{LemM2} will be
\begin{lemma}\label{LemM6} Assume that $\alpha,\beta\in L^2_\per(\Box)$, %$\langle \alpha\beta\rangle=0$,
 $\alpha_\e(x)=\alpha(\e^{-1}x)$, $\beta_\e(x)=\beta(\e^{-1}x)$,
  $\varphi\in L^2(\rd)$, $\psi\in H^1(\rd)$. Then
  \begin{equation}\label{m.14}
|( \alpha_\e S^\e\varphi,\beta_\e S^\e \psi) - \langle \alpha\beta\rangle
( \varphi,\psi)
|\le C\e
\langle \alpha^2\rangle^{1/2}\langle \beta^2\rangle^{1/2}
\|\varphi\|_{L^2(\rd)}\|\nab\psi\|_{L^2(\rd)},\quad C=const(d).
\end{equation}
\end{lemma} 

 The proof of the last three lemmas is given 
  in \S 7.
 \bigskip

\section{Proof of the main results}
We now prove theorems 
\ref{Th2.1} and \ref{Th2.2}.

\textbf{5.1. $H^1$-estimates.}
In  what follows, we  use the notation
\begin{equation}\label{3.1}
u^{,\varepsilon}(x):=S^\e u(x),\quad
 N_\e(x):= N(\f{x}{\e}),\quad 
U^\e(x):=N_\e(x)\cdot\nab u^{,\varepsilon}(x). %S^\e u(x)
 \end{equation}
Then the following estimates %(\ref{190}) and (\ref{1900}) 
hold:  
 \begin{equation}\label{3.2}
\|u^\e-u^{,\varepsilon}-\varepsilon U^\e\|_{H^1(\rd)}\le c\e\|f\|_{L^2(\rd)}, \quad 
c=cost(d,\lambda),
\end{equation} 
  \begin{equation}\label{3.3}
\|u^\e-u-\varepsilon U^\e\|_{H^1(\rd)}\le c\e\|f\|_{L^2(\rd)}, \quad 
c=cost(d,\lambda).
\end{equation}
The latter one is, clearly, equivalent to (\ref{10}).

We give here the proof of
(\ref{3.2})  (the other estimate (\ref{3.3}) is its immediate corollary due to the property (\ref{m.3}) of the operator $S^\e$ and the elliptic estimate (\ref{18})). Further, we % make  
systematically use 
 the estimate (\ref{3.2}) itself and  different elements  in  its proof either.

We begin with necessary calculations:
\begin{equation}\label{3.4}
\ds{
\nab(u^{,\e}+\varepsilon U^\e)=
\nab(u^{,\e}+\e N_\e\cdot\nab u^{,\e})=\left(\nab N_\e^j+e^j\right)\f{\partial u^{,\e}}{\partial x_j}+
\e N_\e^j\nab \f{\partial u^{,\e}}{\partial x_j},
}
\atop\ds{a_\e \nab(u^{,\e}+\varepsilon U^\e)-a^0\nab u^{,\e}=
 g_\e^j\f{\partial u^{,\e}}{\partial x_j}+\e a_\e
  N_\e^j\nab \f{\partial u^{,\e}}{\partial x_j}
}
\end{equation}
(we recall that %as usual, 
summation over repeated  indices is assumed from 1 to $d$),
where 
$$\nab N_\e^j(x):=(\nab_y N^j)(\f{x}{\e}), \quad
g_\e^j(x):=g^j(\f{x}{\e}),$$ 
and the %1-periodic 
vector
 $$
g^j(y):=a(y)\left(\nab N^j(y)+e^j\right)-a^0 e^j,\quad j=1,\ldots,d,
$$
is defined in (\ref{g}). 
 We recall that $ g^j$ is solenoidal and has zero mean value
 (see (\ref{sol})). 
From  (\ref{3.4}), we derive 
\begin{equation}\label{3.7}
\ds{
A_0 u^{,\e}-A_\e (u^{,\e}+\varepsilon U^\e)=
\div\,\left(
a_\e \nab(u^{,\e}+\e U^\e)-a^0\nab u^{,\e}
\right)
=r^\e+\div\, R^\e,
}
\atop\ds{r^\e=
 g_\e^j\cdot\nab \f{\partial u^{,\e}}{\partial x_j},\quad 
 R^\e=\e a_\e
  N_\e^j\nab \f{\partial u^{,\e}}{\partial x_j},
}
\end{equation} 
which enables us to estimate the discrepancy of the  approximation 
 $u^{,\e}+\varepsilon U^\e$ to the equation
 (\ref{1}). Namely,
\[
(A_\e+1)(u^{\e}-u^{,\e}-\varepsilon U^\e)=
(A_\e+1)u^{\e}-(A_\e+1)(u^{,\e}+\varepsilon U^\e)=
\]
\[
(A_0+1)u-(A_\e+1)(u^{,\e}+\varepsilon U^\e)=
\]
\[
(A_0+1)u^{,\e}-(A_\e+1)(u^{,\e}+\varepsilon U^\e)+f-f^{,\e}=
\]
\[
A_0u^{,\e}-A_\e(u^{,\e}+\varepsilon U^\e)-\e U^\e+(f-f^{,\e})\stackrel{(\ref{3.7})}=
\]
\begin{equation}\label{res}%\[
r^\e+\div\, R^\e-\e U^\e+(f-f^{,\e})=:F^\e.
\end{equation} %\]
It is easy to show that
\begin{equation}\label{3.30}%$$
\|F^\e\|_{H^{-1}(\rd)}\le C\e\|f\|_{L^2(\rd)},\quad C=const(d,\lambda),
\end{equation}%$$
using Lemma \ref{LemM3} and estimates (\ref{m.3}), (\ref{18}), if the structure of the functions $r^\e$, $ R^\e$, $U^\e$ (see (\ref{3.1}) and (\ref{3.7})) is taken into account.

To obtain (\ref{3.2}) it remains to apply the following energy inequality
\[
\|z^{\e}\|_{ H^{1}(\rd)}\le c\|F^\e\|_{H^{-1}(\rd)},\quad c=const(\lambda),
\]
to the solution of the equation
\[
z^{\e}\in H^{1}(\rd),\quad (A_\e+1)z^{\e}=F^\e,
\]
where $z^{\e}=u^{\e}-u^{,\e}-\varepsilon U^\e$.

\textbf{5.2. 
$L^2$-estimates.} From (\ref{3.3}), % следует $L^2$-оценка
we have, in particular,
 \[% \begin{equation}\label{3.3}
\|u^\e-u-\varepsilon U^\e\|_{L^2(\rd)}\le c\e\|f\|_{L^2(\rd)}, \quad 
c=cost(d,\lambda),
\]
wherefrom  the $L^2$-estimate (\ref{9}) follows,
since $\|U^\e\|\le c\|f\|$ by properties of smoothing.

On the next step  we %try 
would like to estimate the   $L^2$-norm % норму
 $\|u^\e-u-\varepsilon U^\e\|$ more accurately,
investigating the  $L^2$-form
\[
(u^\e-u-\varepsilon U^\e,h), \quad h\in L^2(\rd).
\]
To this end,  insert $ u^\e-u-\varepsilon U^\e$ as a test function
into the integral identity
for the solution of the adjoint equation 
(\ref{1s}) with an arbitrary function $h\in L^2(\rd)$ on the right-hand side.

We recall % preliminarily 
some facts about the adjoint equation (\ref{1s}) and its solution  $v^\e$. First, the homogenized equation associated with
 (\ref{1s}) is of the form
\begin{equation}\label{3.9}
v\in H^{1}(\rd),\quad (A^*_0+1)v=h;
\end{equation}
second, 
 the approximation in  $H^1$-norm  to $v^\e$
can be chosen as
\begin{equation}\label{3.10}
v^{,\e}(x)+\varepsilon V^\e(x),\quad \mbox{ where } V^\e(x)=
\tilde{N}_\e(x)\cdot\nab v^{,\varepsilon}(x),\quad v^{,\e}(x)=S^\e v(x).
\end{equation}
Here $\tilde{N}$ is the vector composed of of the solutions to the adjoint cell problem % из решений сопряженной задачи на ячейке 
(\ref{cps}).
What is more, the following estimate (that is a counterpart of (\ref{3.3}))  is valid
\begin{equation}\label{3.100}
\|v^\e-v^{,\varepsilon}-\varepsilon V^\e\|_{H^1(\rd)}\le c\e\|h\|_{L^2(\rd)},\quad c=const(d,\lambda).
\end{equation}

Thus, we  write the integral identity for
 $v^\e$ with the test function $ u^\e-u-\varepsilon U^\e$ and make  restructuring  in it:
\[
(u^\e-u-\varepsilon U^\e,h)\stackrel{(\ref{1s})}=(u^\e-u-\varepsilon U^\e,(A^*_\e+1)v^\e)=
\]
\[
((A_\e+1)u^\e-(A_\e+1)(u+\varepsilon U^\e),v^\e)=
%\]
%\[
((A_0+1)u-(A_\e+1)(u+\varepsilon U^\e),v^\e)=
\]
\begin{equation}\label{3.13}\ds{
(A_0u^{,\e}-A_\e(u^{,\e}+\varepsilon U^\e),v^\e)+
(A_0(u-u^{,\e}),v^\e)-(A_\e(u-u^{,\e}),v^\e)-\e( U^\e,v^\e)
=:}\atop\ds{
T_1+T_2-T_3-T_4.}
\end{equation}
Our goal is to estimate the terms $T_i$.

We begin with the simplest term
\[
T_4:=\e( U^\e,v^\e)\stackrel{(\ref{3.1})}=\e( N_\e\cdot\nab u^{,\varepsilon},v^\e)\le \e^2C\langle |N|^2\rangle^{1/2}\|\nab u\|\,\|\nab v^\e\|,
\] 
where the final inequality is due to
Lemma \ref{LemM2} (note that  $\langle N\rangle=0$, see the cell problem (\ref{cp})). Hence, 
in view of (\ref{13s}) и (\ref{18s}), we obtain
\begin{equation}\label{3.15}%\[
T_4\cong 0.
\end{equation}%\]
Here and in the sequel, we  use the notation $\cong$ to denote any equality modulo terms $T$ having the following estimate
\[
|T|\le c\e^2\|f\|\,\|h\|,\quad c=const(d,\lambda);
\]
and such terms $T$ will be called  inessential. 

Next, the term $T_3$ in (\ref{3.13}) admits the following presentation: %and subsequent estimating:
\begin{equation}\label{3.16}%\[
T_3=(u-u^{,\e},A^*_\e v^\e)\stackrel{(\ref{1s})}=(u-u^{,\e},h- v^\e)
\cong 0,
\end{equation}%\]
by (\ref{m.7}) and (\ref{18s})). % and (\ref{13s}).
The similar arguments are applicable to  the term $T_2$. Namely,
\[
T_2:=(A_0(u-u^{,\e}), v^\e)\stackrel{(\ref{3})}=(f-f^{,\e},v^\e)+(u-u^{,\e}, v^\e)\stackrel{(\ref{m.7})}
\cong (f-f^{,\e},v^\e).
\]
We engage now the $H^1$-approximation (\ref{3.10}) and continue our changes:
\[
T_2\cong (f-f^{,\e},v^\e-v^{,\e}-\varepsilon V^\e)+
(f-f^{,\e},v^{,\e}+\varepsilon V^\e)\cong (f-f^{,\e},v^{,\e}+\varepsilon V^\e),
\]
where one term has been dropped,  
because it is  inessential  
in view of the estimates
\[
\|f-f^{,\e}\|_{H^{-1}}\stackrel{(\ref{m.3})}\le C\e\|f\|,\quad
\|v^\e-v^{,\varepsilon}-\varepsilon V^\e\|_{H^1}\stackrel{(\ref{3.100})}\le c\e\|h\|.
\]
Therefore, 
\[
T_2\cong  ((A_0+1)(u-u^{,\e}),v^{,\e}+\varepsilon V^\e)=
\]
\[
(u-u^{,\e},(A^*_0+1)v^{,\e})-\e(f^{,\e}, V^\e)+
\e(f, V^\e)=
\]
\[
(u-u^{,\e},h^{,\e})-\e(f^{,\e}, \tilde{N}_\e\cdot\nab v^{,\varepsilon})+
\e(f, V^\e).
\]
There are  inessential terms in this sum: $(u-u^{,\e},h^{,\e})\cong0$
by (\ref{m.7}), and the next term  is inessential by Lemma \ref{LemM2} (note that  $\langle \tilde{N}\rangle=0$, $f\in \ld$, $\nab v^{,\varepsilon}\in H^1(\mathbb{R}^d)$).
Consequently,
\begin{equation}\label{3.17}%\[
T_2\cong \e(f, V^\e).
\end{equation}

We proceed now to the most difficult term  $T_1$ in (\ref{3.13}). Using the presentation (\ref{3.7}),
we write
\begin{equation}\label{3.18}
T_1:=(A_0u^{,\e}-A_\e(u^{,\e}+\varepsilon U^\e),v^\e)=
%\]
%\[
(g_\e^j\cdot\nab\f{\partial u^{,\e}}{\partial x_j},v^\e)-(\e a_\e
  N_\e^j\nab \f{\partial u^{,\e}}{\partial x_j} ,\nab v^\e)=:
  I+II.
\end{equation}
 Engaging the approximation (\ref{3.10}), we have the sum
\[
I=(g_\e^j\cdot\nab\f{\partial u^{,\e}}{\partial x_j},v^\e-v^{,\e}-\varepsilon V^\e)+
(g_\e^j\cdot\nab\f{\partial u^{,\e}}{\partial x_j},v^{,\e}+\varepsilon V^\e),
\]
where the first summand is inessential due to  Lemma \ref{LemM2} and relations (\ref{sol})$_2$, (\ref{3.100}) and (\ref{18}). Hence, using the fact that   $g_\e^j$ is the solenoidal vector, we obtain
\[
I\cong (g_\e^j\cdot\nab\f{\partial u^{,\e}}{\partial x_j},v^{,\e}+\varepsilon V^\e)=
-
(g_\e^j\f{\partial u^{,\e}}{\partial x_j},\nab(v^{,\e}+\varepsilon V^\e))=
\]
\[
-\left(g_\e^j\f{\partial u^{,\e}}{\partial x_j},
\left(\nab \tilde{N}_\e^k+e^k\right)\f{\partial v^{,\e}}{\partial x_k}+
\e \tilde{N}_\e^k\nab \f{\partial v^{,\e}}{\partial x_k}\right)=
\]
\[
-\left(\left(\nab \tilde{N}_\e^k+e^k\right)\cdot g_\e^j\,\f{\partial u^{,\e}}{\partial x_j},
\f{\partial v^{,\e}}{\partial x_k}\right)-
\e \left(\tilde{N}_\e^kg_\e^j\,\f{\partial u^{,\e}}{\partial x_j},\nab \f{\partial v^{,\e}}{\partial x_k}\right),
\]
where the gradient $\nab(v^{,\e}+\varepsilon V^\e)$ 
has been calculated in the same way as in   (\ref{3.4}).

The periodic vector $(\nab \tilde{N}^k+e^k)\cdot g^j$
has zero mean value. In fact,
\[
\langle g^j\cdot (\nab \tilde{N}^k+e^k)\rangle=
\langle g^j\cdot \nab \tilde{N}^k\rangle+\langle g^j\rangle \cdot e^k=0,
\]
thanks to (\ref{sol}). Thereby, Lemma \ref{LemM5}, combined with the elliptic estimates for the solutions $u$ and $v$ of the homogenized equations, yields
\[
\left(\left(\nab \tilde{N}_\e^k+e^k\right)\cdot g_\e^j\,\f{\partial u^{,\e}}{\partial x_j},
\f{\partial v^{,\e}}{\partial x_k}\right)\cong 0,
\]
and, thus,
\begin{equation}\label{3.19}
I\cong -
\e \left(\tilde{N}_\e^kg_\e^j\,\f{\partial u^{,\e}}{\partial x_j},\nab \f{\partial v^{,\e}}{\partial x_k}\right).
\end{equation}

To estimate the term $II$ in (\ref{3.18}) we write it as the sum
\[
II=-\e( a_\e
  N_\e^j\nab \f{\partial u^{,\e}}{\partial x_j}  ,\nab (v^\e -v^{,\e}-\varepsilon V^\e))-\e( a_\e
  N_\e^j\nab \f{\partial u^{,\e}}{\partial x_j}  ,\nab (v^{,\e}+\varepsilon V^\e)).
\]
Here, the first summand is inessential. To show this, we need only to apply the H\"{o}lder inequality, Lemma
 \ref{LemM1}
and (\ref{3.100}). Next, the calculation of the type of (\ref{3.4}) for the gradient $\nab(v^{,\e}+\varepsilon V^\e)$ is made, after which 
\[
II\cong-\e( a_\e
  N_\e^j\nab \f{\partial u^{,\e}}{\partial x_j}  ,\left(\nab \tilde{N}_\e^k+e^k\right)\f{\partial v^{,\e}}{\partial x_k}+
\e \tilde{N}_\e^k\nab \f{\partial v^{,\e}}{\partial x_k})=
\]
\[
-\e( a_\e
  N_\e^j\nab \f{\partial u^{,\e}}{\partial x_j}  ,\left(\nab \tilde{N}_\e^k+e^k\right)\f{\partial v^{,\e}}{\partial x_k})-
  \e^2( a_\e
  N_\e^j\nab \f{\partial u^{,\e}}{\partial x_j}  ,
 \tilde{N}_\e^k\nab \f{\partial v^{,\e}}{\partial x_k}),
\]
where the last term is inessential due to the H\"{o}lder inequality, Lemma
 \ref{LemM1}  and the elliptic estimates for the solutions $u$ and $v$.
Then
\[
II\cong-\e( a_\e
  N_\e^j\nab \f{\partial u^{,\e}}{\partial x_j}  ,\left(\nab \tilde{N}_\e^k+e^k\right)\f{\partial v^{,\e}}{\partial x_k})=
 -\e( 
  N_\e^j\nab \f{\partial u^{,\e}}{\partial x_j}  ,a^*_\e\left(\nab \tilde{N}_\e^k+e^k\right)\f{\partial v^{,\e}}{\partial x_k}) =
  \]
  \[
 -\e( 
  N_\e^j\nab \f{\partial u^{,\e}}{\partial x_j}  ,\tilde{g}_\e^k\f{\partial v^{,\e}}{\partial x_k}+(a^0)^*\nab v^{,\e}), 
  \]
where we have inserted %introduced 
the vector   $\tilde{g}^k$ (see its definition in (\ref{gs}))
using the equality
\[
a^*(\nab\tilde{N}^k+e^k)=\tilde{g}^k+(a^0)^*e^k.
\]
Note that % по лемме \ref{LemM2}
\[
-\e( 
  N_\e^j\nab \f{\partial u^{,\e}}{\partial x_j}  ,
  (a^0)^*\nab v^{,\e})\cong 0,
\]
by  Lemma \ref{LemM2}, since  $\langle  N^j\rangle=0.$ In conclusion, we obtain
\begin{equation}\label{3.20}
II\cong -
\e( 
  N_\e^j\nab \f{\partial u^{,\e}}{\partial x_j}  ,\tilde{g}_\e^k\f{\partial v^{,\e}}{\partial x_k}).
\end{equation}

From (\ref{3.18})--(\ref{3.20}), we derive
\begin{equation}\label{3.201}\ds{
T_1\cong 
-\e \left(\tilde{N}_\e^kg_\e^j\,\f{\partial u^{,\e}}{\partial x_j},\nab \f{\partial v^{,\e}}{\partial x_k}\right)
 -\e( 
  N_\e^j\nab \f{\partial u^{,\e}}{\partial x_j}  ,\tilde{g}_\e^k\f{\partial v^{,\e}}{\partial x_k})=
 }\atop\ds{
    -\e \left(\tilde{N}_\e^k g_\e^j\,\f{\partial u^{,\e}}{\partial x_j},\nab \f{\partial v^{,\e}}{\partial x_k}\right)
 -\e\left( 
  N_\e^k \tilde{g}_\e^j\cdot\nab \f{\partial u^{,\e}}{\partial x_k}  ,\f{\partial v^{,\e}}{\partial x_j}\right).}
\end{equation}%\]

From now on, our reasoning will be  different in  selfadjoint and nonselfadjoint  cases. We  consider these cases separately, beginning with the first one.

1$^\circ$  %In the case of selfadjoint operators,
Suppose that the matrix $a$ is symmetric. Hence 
$\tilde{N}^k=N^k$, $\tilde{g}^j=g^j$, thereby, the last two forms in (\ref{3.201}) contain the same
vectors
\[
b^{jk}:=N^k g^j= N^k \tilde{g}^j=\tilde{N}^k g^j
\]
such that $b^{jk}\in L^2_\per(\Box)$, since $N^k\in L^\infty_\per(\Box)$ due to the  maximum principle valid in the scalar problem. Subsequent investigation  of the term $T_1$ 
can be based on Lemma \ref{LemM2}. But we avoid using the  maximum principle in order to make our arguments universal and independant of it. We  rely on Lemma \ref{LemM6}.  For the latter it is enough to have $b^{jk}\in L^1_\per(\Box)$ with 
$N^k,g^j\in L^2_\per(\Box)$ which surely holds.
So by Lemma \ref{LemM6},
  \[
  T_1\cong 
   -\e \left(\langle b^{jk}\rangle \f{\partial u}{\partial x_j},\nab \f{\partial v}{\partial x_k}\right)
 -\e\left( 
  \langle b^{jk}\rangle\cdot\nab \f{\partial u}{\partial x_k}  ,\f{\partial v}{\partial x_j}\right)=0.
\]
To explain the final equality to zero, note that
 the last two forms 
  contain the same constant vector %инаковыми 
 $\langle b^{jk}\rangle$ and, besides,  the following equality
\begin{equation}\label{3.200}%\[
\left( \f{\partial \varphi}{\partial x_j},
   \f{\partial^2 \psi}{\partial x_i\partial x_k} \right)
=
-\left( 
   \f{\partial^2 \varphi}{\partial x_i\partial x_k}  ,\f{\partial \psi}{\partial x_j}\right)
   %\quad \forall\varphi,\psi\in \C0, % ^2(\R^d)
    % \varphi,\psi\in H^2(\R^d).
\end{equation}%\]
is valid for  $\varphi,\psi\in H^2(\R^d)$.

Thus,  all the terms $T_i$ in (\ref{3.13}) have been considered.
They are shown to be inessential except for $T_2$ (see (\ref{3.17})).
As a result, the equality 
\begin{equation}\label{3.21}
(u^\e-u-\varepsilon U^\e,h)\cong(f,\e V^\e)
\end{equation}
is proved, where, according to  (\ref{3.1}) and (\ref{3.10}),
\[
U^\e(x)=N_\e(x)\cdot S^\e \nab u(x),\quad V^\e(x)=N_\e(x)\cdot S^\e  \nab v(x).
\]
We 
give the operator form to (\ref{3.21}).
Since
\[
u^\e=(A_\e+1)^{-1}f,\quad u=(A_0+1)^{-1}f,\]
\[
\e U^\e=\e N_\e\cdot S^\e  \nab (A_0+1)^{-1}f=:\e\K_\e f,\]
\[ \e V^\e=\e N_\e\cdot S^\e  \nab (A_0+1)^{-1}h=:\e\K_\e h,
\]
we get
\begin{equation}\label{3.22}
\left((A_\e+1)^{-1}f-(A_0+1)^{-1}f-\e\K_\e f-\e(\K_\e)^* f
, h\right)\cong 0.
\end{equation}
Recalling the convention about the notation  $\cong$ (it is given %immediately 
after (\ref{3.15})), we deduce that
%операторную 
\begin{equation}\label{3.023}\ds{
\|(A_\e+1)^{-1}f-(A_0+1)^{-1}f-\e\K_\e f -\e(\K_\e)^*f  \|
%_{\ld\to \ld}
\le C\ve^2\|f\|, }\atop\ds{
\K_\e= N_\e\cdot S^\e  \nab (A_0+1)^{-1},}
\end{equation}
with the constant $ C=const(d,\lambda)$, 
whence
(\ref{2.21}) immediately follows. 

2$^\circ$ In the nonselfadjoint  case, which means that the matrix $a$ is nonsymmetric,
 the term $T_1$ in (\ref{3.13}) cannot be considered as inessential, thereby, it will contribute to the correcting operator. 
Regarding
 the last two forms in (\ref{3.201}), we see 
  $\e$-periodic vectors
$N_\e^k\tilde{g}_\e^j$ and $\tilde{N}_\e^k g_\e^j$ that 
are distinct.
For the corresponding  1-periodic vectors, we introduce their mean values
\begin{equation}\label{3.240}
c^{jk}=\langle N^k\tilde{g}^j\rangle,\quad
\tilde{c}^{jk}=\langle \tilde{N}^k g^j\rangle.
\end{equation}
By definitions of $\tilde{g}^j$, $g^j$ (see  (\ref{g}), (\ref{gs})), we have 
\begin{equation}\label{3.24}
c^{jk}=\langle N^k a^*(\nab\tilde{N}^j+e^j)\rangle,\quad
\tilde{c}^{jk}=\langle \tilde{N}^k a(\nab{N}^j+e^j)\rangle.
\end{equation}
For instance,  
\[
\tilde{c}^{jk}=
\langle \tilde{N}^k g^j\rangle\stackrel{(\ref{g})}=\langle \tilde{N}^k( a(\nab{N}^j+e^j)-a^0e^j)\rangle=\]
\[\langle \tilde{N}^k a(\nab{N}^j+e^j)\rangle- \langle \tilde{N}^k\rangle a^0e^j\stackrel{(\ref{cps})}
=\langle \tilde{N}^k a(\nab{N}^j+e^j)\rangle.
\]

The same arguments that were used in the selfadjoint case now show
that $\e$-periodic vectors $N_\e^k\tilde{g}_\e^j$  and $\tilde{N}_\e^k g_\e^j$ in (\ref{3.201}) can be replaced with the constant vectors
$c^{jk}$ and
$\tilde{c}^{jk}$, defined in (\ref{3.240}) and (\ref{3.24}), with a negligible error.
As s result,
\begin{equation}\label{3.25}\ds{
T_1\cong 
-\e \left(\tilde{c}^{jk}\,\f{\partial u}{\partial x_j},\nab \f{\partial v}{\partial x_k}\right)
 -\e( 
  c^{jk}\cdot\nab \f{\partial u}{\partial x_j}  ,\f{\partial v}{\partial x_k})=
}\atop\ds
{
\e \left( u,\tilde{c}_i^{jk}\, \f{\partial^3 v}{\partial x_j\partial x_i\partial x_k}\right)
 +\e\left( 
  c_i^{jk} \f{\partial^3 u}{\partial x_j\partial x_i\partial x_k}  , v\right):=
\e (  u,\tilde{L}v )  
 + \e \left(L  u,v \right),}
\end{equation}%\]
where we have introduced the third-order differential operators $L$ and $\tilde{L}$ with the constant coefficients and, thus, completed studying 
the term $T_1$ in (\ref{3.13}).

Gathering the essential terms in (\ref{3.13}),
we obtain
\begin{equation}\label{3.21ns}
(u^\e-u-\varepsilon U^\e,h)\cong(f,\e V^\e)+ \e \left(L  u,v \right)+\e \left(  u,\tilde{L}v \right)  ,
\end{equation}
which should be rewritten in the operator form. To this end, recall that
\[
u^\e=(A_\e+1)^{-1}f,\quad u=(A_0+1)^{-1}f,
\quad v=(A^*_0+1)^{-1}h,\]
%\begin{equation}\label{3.240}%
\[
 U^\e= N_\e\cdot S^\e  \nab (A_0+1)^{-1}f=:\K_\e f,
\quad  V^\e=\tilde{N}_\e\cdot S^\e  \nab (A^*_0+1)^{-1}h=:\tilde{\K}_\e h
%\end{equation}%
\]
and  coin a new operator 
\[
(A_0+1)^{-1}\left(L
+\tilde{L}^*\right) (A_0+1)^{-1}f=:%\tilde
\L f,
\]
where
\begin{equation}\label{3.250}%\[
L
+\tilde{L}^*\stackrel{(\ref{3.25})}=
\left( c_i^{jk}-\tilde{c}_i^{jk}\right)\f{\partial^3 }{\partial x_j\partial x_i\partial x_k}. %:=\mathcal{M}
\end{equation}%\]
Then
\begin{equation}\label{3.260}
\left((A_\e+1)^{-1}f-(A_0+1)^{-1}f-\e\K_\e f-\e(\tilde{\K}_\e)^* f-\e\L f
, h\right)\cong 0.
\end{equation}
 Finally, recalling our convention about the notation
 $\cong$,  
 we establish 
the estimate
\begin{equation}\label{3.23}%\ds{
\|(A_\e+1)^{-1}f-(A_0+1)^{-1}f-\e\K_\e f -\e(\tilde{\K}_\e)^*f -\e \L f \|
%_{\ld\to \ld}
\le C\ve^2\|f\| 
%}\atop\ds{
%\K_\e= N_\e\cdot S^\e  \nab (A_0+1)^{-1},\quad
%\tilde{\K}_\e=\tilde{N}_\e\cdot S^\e  \nab (A^*_0+1)^{-1} %}
\end{equation}
with 
the constant $ C=const(d,\lambda),$ %Из (\ref{3.23}) следует 
whence the estimate (\ref{2.021}) follows.

%\noindent
%\textbf{Замечание.} 

%\textbf{4.3. } 
Since the solutions of the cell problems (\ref{cp}) and 
(\ref{cps}) belong to the space $L^\infty(\Box)$ (recall that we consider the scalar case under the condition (\ref{2})), the functions
$N_\e\cdot  \nab u$ and 
$\tilde{N}_\e\cdot   \nab v$ are well defined as elements of 
$\ld$. 
 If we omit smoothing in the definitions  (\ref{3.1}) and (\ref{3.10}), 
we obtain
$N_\e\cdot  \nab u$ and
$\tilde{N}_\e\cdot   \nab v$ in the place of the correctors
$U^\e$ and $V^\e$. % Moreover, these approximations can be replaced 
   Replacing $U^\e$ и $V^\e$ in (\ref{3.21ns})
with their simplified counterparts $N_\e\cdot  \nab u$ and
$N_\e\cdot   \nab v$, we get  an admissible error, 
 due to the property 
(\ref{m.2}) for the operator $ S^\e$ and the elliptic estimates for $u$ and  $v$  (see (\ref{18}) и (\ref{18s})). Hence we 
 successively find  (\ref{3.260})
and (\ref{3.23}), where smoothing is omitted in  $\K_\e$ and $\tilde{\K}_\e$, which gives (\ref{2.22}). This completes the proof of Theorem \ref{Th2.2}.

\section{Some extension}
 %\setcounter{theorem}{0} \setcounter{equation}{0}
%Diffusion operator with unbounded coefficients
%

\textbf{6.1. Problem setup.} Let us try to weaken
the
conditions (\ref{2}) on the matrix $a(y)$    so that
the main results of  \S 3 (we have in mind the operator
$L^2$-estimates %(\ref{8}),
 (\ref{2.21}) and  (\ref{2.22})) will be still valid.
Assuming that the measurable 1-periodic
matrix $a(y)$ is not symmetric, we decompose it into the symmetric and skew-symmetric
parts:
\beq\label{2s0}%\[
a(y)=a^s(y)+b(y),
\eeq%\]
and we suppose that the symmetric part $a^s$ satisfies the elliptic inequality 
\beq\label{2s}
\lambda |\xi|^2\le {a^s}\xi\cdot\xi
%\quad {a}\xi\cdot\eta
\le\lambda^{-1}|\xi|^2
\quad\forall \xi\in\R^d,\quad\lambda>0.
\eeq

A condition on the skew-symmetric part $b$ is  imposed to ensure, first of all, the unique solvability of the resolvent equation (\ref{1}).
According to the Lax–-Milgram lemma, for this purpose it is sufficient to ensure the boundedness of form
$(a \nab u,\nab \varphi)_{L^2(\R^d)}$  with respect to $u, \varphi \in  H^1 (\R^d)$:
% In other words, we need an estimate
\begin{equation}\label{b.1}
(a\nab u,\nab\varphi)_{L^2(\rd)}\le c_0\|\nab u\|_{L^2(\rd)} \|\nab \varphi\|_{L^2(\rd)}.
\end{equation}
Note that the coercivity of this form, that is, the inequality
\[
(a\nab u,\nab u)_{L^2(\rd)}\ge \lambda\|\nab u\|^2_{L^2(\rd)}
\]
is already ensured 
by the ellipticity of the matrix $a^s$. Moreover, (\ref{2s}) implies also the
boundedness of the $L^2$-form with the matrix $a^s$, and so we need to investigate only  %skew-symmetric
the form
\begin{equation}\label{b.2}%\[
\ds{
(b\nab u,\nab\varphi)_{L^2(\rd)}=
\ild b_{ij}\f{\partial u}{\partial x_i}
\f{\partial \varphi}{\partial x_j}\,dx
=\f{1}{2}\ild b_{ij}\left(\f{\partial u}{\partial x_i}
\f{\partial \varphi}{\partial x_j}-\f{\partial u}{\partial x_j}
\f{\partial \varphi}{\partial x_i}\right)\,dx=}
\atop\ds{
\f{1}{2}\ild  b_{ij} I_{ij}(u,\varphi)\,dx.}
\end{equation}%\]

Note that the necessary and sufficient 
conditions on  the matrix $a$ for the continuity property (\ref{b.1}) were investigated in \cite{MV}. Dealing with homogenization, we have to reproduce some details of this investigation.

The skew-symmetric difference $I_{ij}(u,\varphi)$ in parentheses
of (\ref{b.2}) have “better than expected” regularity: it 
 belongs  surely to the space
$L^1(\rd)$, but 
%the ‘compensated’ regularity comes into play: 
the algebraic structure
makes this non-linear expression lie in the narrower Hardy space
  \[
\mathcal{H}^1(\rd)=\{f\in L^1(\rd):\,R_jf\in  L^1(\rd),\, 1\le j\le d\},
\]
where 
 $R_j{=}\f{\partial }{\partial x_j}(-\Delta)^{-1/2}$
 are the Riesz operators  (see Proposition 4.4 in  \cite{UMN} which is proved relying upon the results from \cite{CLMS}).
 
 The dual of the space $\mathcal{H}^1(\rd)$ is the space $BMO$ (bounded mean oscillation) \cite{St}. We recall that a measurable
function $g$ on $\rd$ lies in $BMO$ if
\begin{equation}\label{bmo}%\[
\|g\|_{BMO}=\sup \sr B|g-g_B|\,dx\le \infty\quad\mbox{ with }\quad g_B{=}\sr B g\,dx{=}\f{1}{|B|}\int_B g\,dx,
\end{equation}%\]
where the supremum is taken over all balls $B{\subset}\rd$.
% and $g_B{=}\sr B g\,dx{=}\f{1}{|B|}\int_B g\,dx$.
An equivalent definition is obtained if we replace the balls by cubes.
Obviously, elements of the space $BMO$ are defined up to a constant.

By duality arguments (see more details in \cite{UMN}, \S 4), we come to 
\begin{prop}\label{PrB}
If the entries of the matrix $b$ lie in $BMO$, then the form (\ref{b.2})
is bounded:
\begin{equation}\label{bform}%     \[
     (b\nab u,\nab\varphi)_{L^2(\rd)}\le c_0 \|\nab u\|_{L^2(\rd)}
\|\nab \varphi\|_{L^2(\rd)},
\end{equation}%
where the constant $c_0$ depends only on the norm $\|b\|_{BMO}$.
 \end{prop}
 
Thus, from now on we assume:

(\textbf{C}) \textit{ the symmetric part of the matrix $a$ satisfies the ellipticity condition (\ref{2s});\\
 $\quad{}\quad{} \quad{}\quad\,\,$ its skew-symmetric part $b$ 
belongs to the space} $BMO$.\\
 Then 
the whole form $(a\nab u,\nab\varphi)_{L^2(\rd)}$ is bounded, and 
the estimate (\ref{b.1}) holds with the constant
$c_0$ depending only on $\lambda$ and $ \|b\|_{BMO}$. 
%We note that 
A homothety does not change the $BMO$ norm:
if $b_\varepsilon(x) = b(x/\varepsilon)$, then $\|b_\e\|_{BMO} =\|b\|_{BMO}$. Hence, the form with an $\e$-periodic
matrix $a_\e(x)$, namely, $(a_\e\nab u,\nab\varphi)_{L^2(\rd)}$, is bounded  and satisfies an estimate of type
(\ref{b.1})  with the same constant $c_0$.

Therefore, Equation (\ref{1}) is uniquely solvable and the uniform (in $\e$)  estimate of the type (\ref{5}) is valid for its solution.
 Parallelly, one can show that the cell problem (\ref{cp}) (and also (\ref{cps})) is well posed, thereby, the homogenized matrix $a_0$ is well defined in (\ref{13}) in terms of the solutions $ N^j$ to (\ref{cp}) (see details in \cite{UMN}, \S 4).

\textbf{6.2. Estimates of order $\e$.}  Under condition (\textbf{C}), all the homogenization results stated in \S 2 remain %are 
 true, including the estimate (\ref{11}). The maximum principle
holds for the cell problems (\ref{cp}), (\ref{cps}), and its solutions $N^j$ and $\tilde{N}^k$ belong to $L^\infty(\Box)$ (for proof see arguments in \cite{SSSZ})).

 To justify the operator estimates  (\ref{8}) and (\ref{10}) under condition (\textbf{C}), look through the reasoning in \S 5, where we
derive the main estimate  (\ref{3.2})
 from which (\ref{8}) and (\ref{10}) easily follow.  One of the key points here is the estimate (\ref{3.30}) for the residual $F^\e$ defined in (\ref{res}). 
 To obtain this estimate  we are to benefit from  Lemma 
 \ref{LemM3}. For this purpose, the terms $r^\e$ and $R^\e$ defined in (\ref{3.7}) should have a proper structure, that is, its components  
 $g^j$ and $aN^j$ should be sufficiently regular, namely,
 \begin{equation}\label{re}
 g^j,\, aN^j \in L^2(\Box).
 \end{equation}
 
By definition (\ref{g}), to show $ g^j \in L^2(\Box)$
%the first inclusion in (\ref{re})
 we  are to invoke
%is a corollary of
 \begin{prop}\label{LemB1}
If $N^j$ is the solution of the problem (\ref{cp}), then $a\nab N^j$ belongs to
$L^2(\Box)$ and satisfies the estimate $\|{a}\nab N^j\|_{L^2(\Box)}{\le}C$, where the constant $C$ depends only on $\lambda$ and $ \|b\|_{BMO}$.
%has the same type as above.
\end{prop}
This assertion  is proved  in \cite{UMN}  relying on the higher integrability of the gradient $\nab N^j$, that is,   $\nab N^j\in L^{2+\delta}_\per(\Box)$ for some $\delta>0$ (see Lemma 4.2 in \cite{UMN}), and  the John--Nirenberg inequality
\begin{equation}\label{hiin}%\[
\sr B|g-g_B|^p\,dx\le c_p \|g\|^p_{BMO}\quad\forall p>1,
\end{equation}%\]
which stems from (\ref{bmo}). We apply (\ref{hiin}) to the 
 function $b$ (the skew-symmetric part of the matrix $a$) on the unit cube $B=\Box$.
 It is appropriate here to refer to the fact that
the form (\ref{b.2}) will not change its value on subtracting a constant skew-symmetric
matrix $C$ from $b$. In the case of $b \in BMO$, a suitable integral mean $g_B$ (see (\ref{bmo}))
is %frequently 
taken for this constant matrix, which allows one to invoke the 
John--Nirenberg inequality.

%The second inclusion in (\ref{re}) is proved 
To show $aN^j \in L^2(\Box)$ it is enough to apply the assertion
%a corollary of 
(\ref{hiin}) with respect to   the 
 matrix $b$ and the higher integrability of $N^j$ by the Sobolev embedding theorem (or the deeper property  $N^j \in L^\infty(\Box)$).
 
 In conclusion of this Subsection, note that the more detailed proof of the estimates (\ref{8}) and (\ref{10}) for the operator $A_\e$  with the coefficients from $BMO$ is given in  \cite{UMN}.

\textbf{6.3. $L^2$-estimate of order $\e^2$.} Assuming the condition (\textbf{C}) on the matrix $a$ stated in Subsection 6.1, 
let us %prove
show that  the operator estimate (\ref{2.21}) remains valid.
%for the resolvent of the equation  (\ref{1}).
 We can repeat without any changes  reasoning of \S6 up to the “equality” (\ref{3.201}), in particular, taking into account  (\ref{re}),        which 
%is not evident, in the case of the unbounded matrix  $a$, 
 has been already explained (see  the preceding  subsection). Then, following the lines of the nonselfadjoint case, we come to (\ref{3.25}) by Lemma  \ref{LemM6} and afterwards  duplicate the end of the proof of the estimate (\ref{2.21}) in \S6.

We formulate finally the main result of this section that has been just verified.
 \begin{teor}\label{Th6.1} Let the matrix $a$ in (\ref{1}) satisfy
 (\ref{2s0}), (\ref{2s}) with the skew-symmetric part $b\in BMO$. 
 Then there holds the    estimate (\ref{2.021}) with the correcting operators defined in (\ref{2.121}) and the constant $C$ in the right-hand side, which  depends only on the dimension $d$, the ellipticity constant $ \lambda$ in (\ref{2s}) and the norm $\|b\|_{BMO}$.
\end{teor} 

 \section{Auxiliaries}
In this section  we give the proof of some properties of the Steklov smoothing operator $S^\e$ formulated in \S4.

\textbf{Proof of Lemma \ref{LemM4}.}
To estimate the form
$
I:=( b_\e S^\e\varphi,S^\e\psi),
$
 where
\[
S^\e\varphi(x)=\ilb\varphi(x-\e\omega)\,d\omega,\quad  S^\e\psi(x)=\ilb\psi(x-\e\sigma)\,d\sigma,
\]
 we make standard transformations: 
\[
I=\ild\ilb\ilb b(\f{x}{\varepsilon})\varphi(x-\varepsilon\omega)
\psi(x-\varepsilon\sigma)\,d\omega\,d\sigma\,dx=
\]
\[\ild\ilb\ilb b(\f{x}{\varepsilon}+\omega+\sigma)\varphi(x+\e\sigma)\psi(x+\varepsilon\omega)\,d\omega\,d\sigma\,dx=\]
\[
\ild\ilb\ilb b(\f{x}{\varepsilon}+\omega+\sigma)\varphi(x+\e\sigma)\left(\psi(x+\varepsilon\omega)-
\psi(x)\right)
\,d\omega \,d\sigma\,dx=
\]
\[
\ild\ilb\ilb b(\f{x}{\varepsilon}+\omega+\sigma)\varphi(x+\e\sigma)
\left(\int\limits_{0}^1\nab\psi(x+t\varepsilon\omega)\cdot 
\varepsilon\omega\,dt \right)
\,d\omega\,d\sigma\,dx
\]
where we have used the condition
$\langle b\rangle=0$ and the integral representation
\begin{equation}\label{m.10}%\[
\psi(x+h)-\psi(x)=
\int\limits_{0}^1\nab\psi(x+th)\cdot h\,dt.
\end{equation}
We continue the  standard transformations: 
\[
I=
\int\limits_{0}^1\ild\ilb\ilb b(\f{x}{\varepsilon}+\omega+\sigma)(\varphi(x+\varepsilon\sigma)-\varphi(x))\nab\psi(x+t\varepsilon\omega)\cdot 
\varepsilon\omega 
\,d\omega\,d\sigma\,dx\,dt=
\]
\[
\int\limits_{0}^1\ild\ilb\ilb b(\f{x}{\varepsilon}+\omega+\sigma)\left(
\int\limits_{0}^1\nab\varphi(x+s\varepsilon\sigma)\cdot \varepsilon\sigma\,ds
\right)\nab\psi(x+t\varepsilon\omega)\cdot 
\varepsilon\omega 
\,d\omega\,d\sigma\,dx\,dt=
\]
\[
\varepsilon^2\int\limits_{0}^1\int\limits_{0}^1\ild\ilb\ilb b(\f{x}{\varepsilon}+\omega+\sigma)\left(
\nab\varphi(x+s\varepsilon\sigma)\cdot\sigma\right)
\nab\psi(x+t\varepsilon\omega)\cdot 
\omega 
\,d\omega\,d\sigma\,dx\,dt\,ds,
\]
where we again use the condition
$\langle b\rangle=0$ and an integral representation
 for the difference $\varphi(x+\varepsilon\sigma)-\varphi(x)$ similar to (\ref{m.10}).
Applying the  H\"{o}lder inequality to the last multi-dimensional integral, we find
 \begin{equation}\label{m.11}\ds{
I^2\le
\varepsilon^4\int\limits_{0}^1\int\limits_{0}^1\ild\ilb\ilb |b(\f{x}{\varepsilon}+\omega+\sigma)|^2|
\nab\varphi(x+s\varepsilon\sigma)\cdot\sigma|^2\,d\omega\,d\sigma\,dx\,dt\,ds\times
}
\atop\ds{
\int\limits_{0}^1\int\limits_{0}^1\ild\ilb\ilb |\nab\psi(x+t\varepsilon\omega)\cdot 
\omega |^2\,d\omega\,d\sigma\,dx\,dt\,ds,
}
\end{equation}
where both integral factors can be easily estimated:
 \begin{equation}\label{m.12}%\[
I^2\le
\varepsilon^4C\langle |b|^2\rangle \|\nab\varphi\|^2_{L^2(\rd)} \|\nab\psi\|^2_{L^2(\rd)},\quad C=const(d).
\end{equation}%\]
Hence we obtain the estimate (\ref{m.8}). The lemma is proved.
 
%%%%%%%%%%%%%%
\textbf{Proof of Lemma \ref{LemM5}.}
Deriving the estimate  (\ref{m.13}), one can assume that
$\varphi,\psi\in \C0$ and, considering the 
oscillating factor $b=\alpha\beta$, repeat the  standard transformations  of the form $I$ made above up to formula  (\ref{m.11}). Before we use the  H\"{o}lder inequality, we recall that $b=\alpha\beta$ and distribute the functions  $\alpha$ и  $\beta$ among the different integral factors.
Thus, instead of (\ref{m.11}), we come to the inequality
 % \begin{equation}\label{m.11}
 $$\ds{
I^2\le
\varepsilon^4\int\limits_{0}^1\int\limits_{0}^1\ild\ilb\ilb |\alpha(\f{x}{\varepsilon}+\omega+\sigma)|^2|
\nab\varphi(x+s\varepsilon\sigma)\cdot\sigma|^2\,d\omega\,d\sigma\,dx\,dt\,ds\times
}
\atop\ds{
\int\limits_{0}^1\int\limits_{0}^1\ild\ilb\ilb 
|\beta(\f{x}{\varepsilon}+\omega+\sigma)|^2|\nab\psi(x+t\varepsilon\omega)\cdot 
\omega |^2\,d\omega\,d\sigma\,dx\,dt\,ds.
}
$$%\end{equation}
Here, both integral factors can be easily estimated and, instead of (\ref{m.12}), we obtain
\[
I^2\le
\varepsilon^4C\langle |\alpha|^2\rangle \|\nab\varphi\|^2_{L^2(\rd)} \langle |\beta|^2\rangle  \|\nab\psi\|^2_{L^2(\rd)},\quad C=const(d),
\]
which is equivalent to (\ref{m.13}). The lemma is proved.

\textbf{Proof of Lemma \ref{LemM6}.} The same transformations that are used in the proof of Lemma \ref{LemM4} yield
\[
I:=( \alpha_\e S^\e\varphi,\beta_\e S^\e \psi)=\ild\ilb\ilb b(\f{x}{\varepsilon})\varphi(x-\varepsilon\omega)
\psi(x-\varepsilon\sigma)\,d\omega\,d\sigma\,dx=
\]
\[\ild\ilb\ilb b(\f{x}{\varepsilon}+\omega+\sigma)\varphi(x+\e\sigma)\psi(x+\varepsilon\omega)\,d\omega\,d\sigma\,dx,\]
where we set $b=\alpha\beta$. Decomposing $\psi(x+\varepsilon\omega)=
\psi(x)+(\psi(x+\varepsilon\omega)-\psi(x))$, we write the representation
$%\[
I=I_1+I_2
$, %\] 
where
\[
I_1:=\ild\ilb\ilb b(\f{x}{\varepsilon}+\omega+\sigma)\varphi(x+\e\sigma)\psi(x)\,d\omega\,d\sigma\,dx=\langle b\rangle
\ild\ilb \varphi(x+\e\sigma)\psi(x)\,d\sigma\,dx=
\]
\[
\langle b\rangle
\ild\ilb \varphi(x)\psi(x-\e\sigma)\,d\sigma\,dx=\langle b\rangle(\varphi,S^\e\psi)
\]
and 
\[
I_2:=\ild\ilb\ilb b(\f{x}{\varepsilon}+\omega+\sigma)\varphi(x+\e\sigma)(\psi(x+\varepsilon\omega)-\psi(x))\,d\omega\,d\sigma\,dx=
\]
\[
\ild\ilb\ilb b(\f{x}{\varepsilon}+\omega+\sigma)\varphi(x+\e\sigma)\int\limits_{0}^1\nab\psi(x+t\varepsilon\omega)\cdot 
\varepsilon\omega\,dt 
\,d\omega\,d\sigma\,dx
\]
%after using 
thanks to the integral formula (\ref{m.10}). By arguments used  in the proof of Lemma \ref{LemM5}, we show that
\[
I_2\le C\e
\langle \alpha^2\rangle^{1/2}\langle \beta^2\rangle^{1/2}
\|\varphi\|_{L^2(\rd)}\|\nab\psi\|_{L^2(\rd)}.
\]

Returning to $I_1$, it is clear that
$%\[
I_1=\langle b\rangle(\varphi,\psi)+\langle b\rangle(\varphi,S^\e\psi-\psi)
$, %\]
where, by properties of smoothing, the second summand admits the estimate from above with the same majorant as  in (\ref{m.14}).
Eventually, gathering all the relations proved above, we come to (\ref{m.14}). The lemma is proved.


\begin{thebibliography}{99}

%%%1
 \bibitem{BLP}  A. Bensoussan,  J.L. Lions,  G. Papanicolaou, \emph{Asymptotic  Analysis for
Periodic Structures,} North Holland, Amsterdam (1978).

%%%2
\bibitem{SP} E. Sanchez-Palencia, \emph{Non-homogeneous media and vibration theory,} Lecture Notes
in Phys., vol. 127, Springer-Verlag, Berlin–New York 1980.

%%%3
\bibitem {BP} 
N. Bakhvalov and G. Panasenko, Homogenisation: Averaging Processes in Periodic
Media: Mathematical Problems in the Mechanics of Composite Materials, Nauka,
Moscow, 1984 (in Russian); Kluwer Academic, Dordrecht, 1989 (in English).

%%4
\bibitem{ZKO} 
V. V. Jikov, S. M. Kozlov and O. A. Oleinik, Homogenization of Differential
Operators and Integral Functionals, Springer, Berlin, 1994.

%%%%5
\bibitem{BS}  M.\,Sh. Birman and T.\,A. Suslina,
"Second order periodic differential operators. Threshold properties and homogenization"\!,  \textit{Algebra i Analiz}, {\bf 15}(5), 1–-108   (2003)(in Russian);  \textit{St. Petersburg
Math. J.,} {\bf 15}, 639–-714  (2004) (in English).

%%%6
\bibitem{Zh1} V.\,V. Zhikov, "On operator estimates in homogenization theory"\!,  \textit{Dokl. Acad. Nauk,}
{\bf 403}(3),  305–-308  (2005) (in Russian);  \textit{Dokl. Math.}, {\bf 72}(1),  535–-538  (2005) (in English).


%%%7
\bibitem{Zh05}	V.\,V. Zhikov, "Spectral method in homogenization theory"\!,  \textit{Proc. Steklov Inst. Math.}, {\bf 250},
85–-94 (2005).

%%%%8
\bibitem{BS05} 
M.\,Sh. Birman and T.\,A. Suslina, "Homogenization with corrector term for periodic elliptic differential operators\!,  \textit{Algebra i Analiz,} {\bf 17}(6), 1–-104  (2005) (in Russian);  \textit{St. Petersbg. Math. J.}, {\bf 17}(6), 897–-973 (2006) (in English).

%%%%5
\bibitem{ZhP05}  V.\,V. Zhikov,  S.\,E. Pastukhova,  "On operator estimates for some problems in homogenization theory"\!,  
 \textit{ Russian Journal of Math.Physics},  {\bf 12}, №4,  515--524 (2005).

\bibitem{ZhP06}  V.\,V. Zhikov,  S.\,E. Pastukhova,  "Estimates of homogenization for a parabolic equation with periodic coefficients"\!,  
 \textit{ Russian Journal of Math.Physics},  {\bf 13}, №4,  224--237 (2006).

%%%%%%%
 \bibitem{Zh2} V. V. Zhikov, “Some estimates from homogenization theory,"\!, \textit{ Dokl. Acad. Nauk,} {\bf 406}(5), (2006),  597--601 (in Russian);  \textit{Dokl. Math.}, {\bf 73}(1), 96–-99
(2006)  (in English).

 
 \bibitem{Pa06} S.\,E. Pastukhova, "Some estimates from homogenized elasticity problems"\!,  \textit{Dokl. Acad. Nauk,}  {\bf 406}(5) , 604–-608 (2006) (in Russian); 
  \textit{Dokl. Math.},  {\bf 73}(1),
102–-106   (2006) (in English).
 

%%%%%%13 
 \bibitem{CPZ} G. Cardone,   S.\,E. Pastukhova and V.\,V. Zhikov, "Some estimates for nonlinear homogenization"\!, \textit{ Rend. Accad. Naz. Sci. XL Mem. Mat. Appl.},  {\bf 29},   101--110 (2005).
 
  \bibitem{ZPT} V.\,V. Zhikov, S.\,E. Pastukhova and S.\,V. Tikhomirova, "On the homogenization of degenerate
elliptic equations"\!,  \textit{Dokl. Acad. Nauk,} {\bf 410}(5), 587--591 (2006)
(in Russian); 
 \textit{ Dokl. Math.,} {\bf  74}(2), 716--720 (2006)  (in English).

 
 \bibitem{ZP08} V.\,V. Zhikov and S.\,E. Pastukhova, "Homogenization of degenerate elliptic
equations"\!,\textit{ Sib. Math. J.}, {\bf 49}(1), 101–-124  (2008) (in Russian); \textit{ Sib. Math. J.},  {\bf 49}(1), 80–-101  (2008) (in English).

 \bibitem{PT} S.\,E. Pastukhova and R.\,N. Tikhomirov, "Operator estimates in reiterated and locally
periodic homogenization"\!, \textit{Dokl. Acad. Nauk,} {\bf 415}(3),  304–-309 (2007) (in Russian);
 \textit{Dokl. Math.,} {\bf 76}, 548–-553 (2007) (in English).  


 \bibitem{PaSt} S.\,E. Pastukhova, "Operator estimates in nonlinear problems of reiterated homogenization"\!,  
 \textit{Tr. Mat. Inst. Steklova,}  {\bf 261},  220–-233  (2008) (in Russian);  \textit{Proc. Steklov Inst. Math.,} {\bf 261}, 214–228  (2008) (in English). 

 \bibitem{PAsAn} S.\,E. Pastukhova, "Estimates in homogenization of parabolic equations with locally periodic coefficients"\!,  
 \textit{  Asymptot. Anal}, {\bf 66}(3-4), 207–-228 (2010).

 \bibitem{PaFAA}  S.\,E. Pastukhova, "Approximation of the exponential of
a diffusion operator with multiscale coefficients"\!,  
 \textit{Funktsionalnyi Analiz i Ego Prilozheniya,} {\bf 48}(3),
34–-51   (2014) (in Russian); \textit{ Funct. Anal. Appl.}, {\bf 48}(3),
183–-197 (2014)  (in English).


 
\bibitem{ZhP15}
V.\,V. Zhikov  and  S.\,E. Pastukhova,  "Homogenization estimates of operator type for an elliptic equation with quasiperiodic coefficients"\!,  
 \textit{ Russian Journal of Math.Physics},  {\bf 22}(4), 
  264–-278 (2015). 
  
  %%%%%
\bibitem{PAA16}
S.\,E. Pastukhova, "Operator error estimates for homogenization of fourth order elliptic equations"\!, \textit{ Algebra i Analiz}, {\bf 28}(2), 
 204–-226 (2016) (in Russian);  \textit{St. Petersburg Math. J.,} {\bf 28}(2), 273–-289  (2017) (in English).

\bibitem{P16}
S.\,E. Pastukhova, 
 "Estimates in homogenization of higher-order elliptic operators"\!,
    \textit{    Applicable Analysis} \textbf{95}(7-9), 1449--1466 (2016).
  
\bibitem{UMN} V.\,V. Zhikov  and S.\,E. Pastukhova, "Operator estimates in homogenization theory"\!,  
 \textit{ Uspekhi Mat. Nauk,} {\bf 71}(3),  27–-122  (2016) (in Russian);  \textit{Russian Math. Surveys,} {\bf 71},  417–-511  (2016) (in English).
   


\bibitem{P17} S.\,E. Pastukhova, "Operator Estimates in Homogenization of Elliptic Systems of Equations"\!,  
 \textit{Journal of Math. Sciences},
 % {\bf 89},  99--112 (2017).
{\bf 226}(4),  445--461  (2017).





 
\bibitem{PT17} S.\,E. Pastukhova and R.\,N. Tikhomirov, "Operator-type estimates in homogenization of
elliptic equations with lower order terms"\!, \textit{ Algebra i Analiz}, {\bf 29}(5),  179--207 (2017) (in Russian);
 \textit{ St. Petersbg. Math. J.}, {\bf 29}(5), 841–-861 (2018)  (in English).




 
\bibitem{P20}   S.\,E. Pastukhova,  "$L^2$-estimates for homogenization of elliptic operators"\!,  
 \textit{ Journal of Math. Sciences},  {\bf 244}(4),  671--685 (2020).



%%%%%%%%%
\bibitem{ZKOH}
V.\,V. Zhikov et al.
"Averaging and G-convergence of differential operators"\!,
 \textit{ Russian Math. Surveys},  {\bf 34}(5),  69–-147, (1979).
% 34:5 , 69–147.

\bibitem{Se1}
N.\,N. Senik, "Homogenization for non-self-adjoint periodic elliptic operators on an
infinite cylinder"\!,
 \textit{ SIAM J. Math. Anal.}, {\bf 49},  874–898  (2017).

\bibitem{Se} N.\,N.  Senik, "Homogenization for non-self-adjoint locally periodic elliptic operators"\!, 1--20 (2017)
 arXiv:1703.02023v2 [math.AP].

\bibitem{MV} V.\,G. Maz’ya and I.\,E. Verbitsky, "Form boundedness of the general second-order
differential operator"\!, \textit{  Comm. Pure Appl. Math.},
 {\bf 59}(9),  1286–-1329  (2006).
  
 
\bibitem{CLMS} R. Coifman, P.\,L. Lions, Y. Meyer, and S. Semmes, "Compensated compactness
and Hardy spaces"\!, \textit{ J. Math. Pures Appl.},
 {\bf 72}(3),  247–-286  (1993).
 %(9) 72:3 (1993), 247–286.

\bibitem{St} E.\,M. Stein, Harmonic analysis: real-variable methods, orthogonality, and oscillatory
integrals, Princeton Math. Ser., vol. 43, Princeton Univ. Press, Princeton, NJ 1993.
%xiv+695 pp.

\bibitem{SSSZ} G. Seregin, L. Silvestre, V. Sverak and  A. Zlatos, "On divergence-free drifts"\!, \textit{ J. Differential
Equations},  {\bf 252}(1),  505--540  (2012).


%\bibitem{ZhFA} V.\,V. Zhikov, "Remarks
%on the uniqueness of a solution of the Dirichlet problem for second-order elliptic
%equations with lower-order terms"\!, \textit{  Funct. Anal. Appl.}, {\bf 38}(3),  173–-183 (2004).
 %38:3 (2004), 173–183.

 

\end{thebibliography}
\end{document}